\documentclass[11pt, reqno, a4paper]{amsart}

\usepackage[normalem]{ulem}
\usepackage[utf8]{inputenc}
\usepackage{a4wide}

\usepackage{algorithm}
\usepackage{algpseudocode}
\usepackage{todonotes}
\usepackage{bbm}

\usepackage{dsfont}

\usepackage{amssymb, amsthm, amsmath, enumerate, physics, mathrsfs, amsfonts}
\usepackage[colorlinks=true, linkcolor=black, citecolor=black, urlcolor=black]{hyperref}
 \usepackage{accents}
\usepackage{xcolor}
\usepackage{fullpage}
\usepackage{graphicx}
\usepackage{pgfplots}
\pgfplotsset{compat=newest}
\usetikzlibrary{pgfplots.groupplots}
\usepackage{array}
\usepackage{subcaption}
\usepackage[edges]{forest}
\usepackage[normalem]{ulem}

\usepackage[numbers]{natbib}

\allowdisplaybreaks
\numberwithin{equation}{section}

\newcommand{\R}{\mathbb{R}}

\newcommand{\C}{\mathcal{C}}

\newcommand{\X}{\mathcal{X}}
\newcommand{\Y}{\mathcal{Y}}
\newcommand{\Pn}{\mathscr{P}}

\DeclareMathOperator*{\argmin}{arg\,min}

\theoremstyle{plain}

\newtheorem{rmk}{Remark}[section]
\theoremstyle{definition}

\newcommand{\response}[1]{\textcolor{black}{#1}}

\begin{document}

\title{A Semi-Discrete Optimal Transport Scheme for the 3D Incompressible Semi-Geostrophic Equations}
\author{Th\'eo Lavier}

\begin{abstract}

We describe a mesh-free three-dimensional numerical scheme  for solving the incompressible semi-geostrophic equations based on semi-discrete optimal transport techniques.
These results generalise previous two-dimensional implementations.
The optimal transport methods we adopt are known for their structural preservation and energy conservation qualities and achieve an excellent level of efficiency and numerical energy-conservation.
We use this scheme to generate numerical simulations of an important cyclone benchmark problem.
To our knowledge, this is the first fully three-dimensional simulation of the semi-geostrophic equations, evidencing semi-discrete optimal transport as a novel, robust numerical tool for meteorological and oceanographic modelling.

\end{abstract}

\maketitle

\section{Introduction}

In this paper we describe what is, to our knowledge, the first mesh-free three-dimensional (3D) numerical scheme providing simulations of incompressible semi-geostrophic (SG) atmospheric flows. 
In particular, we use this scheme to simulate the evolution of an isolated large-scale tropical cyclone, supporting the applicability of the SG equations for modelling atmospheric and oceanic phenomenon.

The SG system is a second-order accurate reduction of the Euler equations valid for modelling large-scale atmospheric flows. 
Its significance in meteorology stems from the fact that the system models the formation of fronts - mathematically, that it has a natural way to admit solutions that continue past singularity formation. 
The recent success and interest in this system is a consequence of its reformulation, due to Brenier and Benamou, as a coupled optimal transport problem. 
Indeed, it is this reformulation that we exploit to devise an energy-conserving approximation that models accurately the formation of fronts and cyclones.

The scheme we present is based on semi-discrete optimal transport. 
\response{
Two-dimensional (2D) reductions of SG dynamics, and of its incompressible Boussinesq parent model, have been approached through several complementary analytical and numerical paradigms. 
In the Lagrangian setting, Benamou et al. and Carlier et al. solve the SG equations with a fully discrete optimal-transport solver regularised by entropy \cite{Benamou:2024,Carlier:2024}.
A counterpart is provided by semi-discrete optimal transport, which supplies a rigorous framework for Cullen’s pioneering \emph{geometric method} and has been investigated analytically both in 2D and 3D and simulated in the canonical 2D Eady slice \cite{Bourne:2022,Egan:2022,Cullen:1984}.
In contrast to the mesh-free optimal-transport approaches, Yamazaki et al. examine the incompressible Eady–Boussinesq slice reproducing quasi-periodic frontogenetic life-cycles using high-resolution Eulerian compatible-finite-element solver that serve as a benchmark for the SG limit \cite{Yamazaki:2017}.
Together, the schemes of \cite{Benamou:2024} and \cite{Yamazaki:2017} provide reference diagnostics and convergence results against which we later qualitatively evaluate our 2D benchmark and 3D simulations.
Although neither their solvers nor the present one are strictly energy-conserving, all retain dissipation that is small relative to the kinetic and potential-energy budgets, so all support credible long-time diagnostics. 
Our semi-discrete optimal-transport framework should therefore be viewed as a complementary alternative for long-time diagnostics.
Sch\"ar and Wernli \cite{Schar:1993}, who originally proposed the initial conditions we investigate, used Eulerian spectral grid with centred finite-difference time stepping to solve the SG equations. 
Fully discrete, entropy-regularised optimal-transport schemes are efficient, energy-stable, and convergent. 
They represent both the source measure and the evolving target measure on a finite discrete support and hence satisfy the semi-geostrophic equations only in a discrete sense. 
One recovers a weak solution of the continuous system only in the joint limit of vanishing entropic parameter \cite{Benamou:2024, Carlier:2024}.
By contrast, in the semi-discrete approach the source measure remains continuous while the target measure is discrete, so the resulting particle trajectories are exact weak solutions of the Lagrangian form of the SG equations for every particle number $N$ \cite{Bourne:2022,Egan:2022}.
} 
Our contribution is to present what is, to the authors' knowledge, the first full 3D simulation of the formation of an isolated cyclone, developing from the benchmark set of initial conditions given in \cite{Schar:1993}.
Our work to implement fully 3D computations is a substantial extension of the method presented in \cite{Egan:2022}. 

In addition to the novel spatial particle discretisation, achieved using semi-discrete optimal transport, we briefly investigate various explicit numerical methods to solve the temporal evolution, focusing on balancing runtime with the relative error in the energy conservation. 
While straight forward, this exploration is novel in terms of integrating an ordinary differential equation (ODE) solver with an optimal transport solver.

We then implement our schemes, first demonstrating that the observed rate of convergence aligns with the theoretical predictions for the space and time discretisation (see \S~\ref{sec:experiments}). 
We then provide insights into large-scale SG flows and how they model the formation and development of atmospheric fronts and cyclones. 

\subsection{Background and motivation}

A central theme in the study of atmospheric and oceanic dynamics is the quest for  models that are both mathematically and numerically tractable and that approximate accurately the fluid motion, at least within a set of specific physical constraints. 
An important system of equations satisfying these requirements is the SG system.
This system, which is derived under the assumptions of hydrostatic and geostrophic balance, is a second-order accurate reduction of the full Euler system and is recognised for its effectiveness in modelling shallow, rotationally-dominated flows characterised by small Rossby numbers. 
This contrasts with the quasi-geostrophic limit which simplifies the dynamics by neglecting ageostrophic terms beyond the first order. 
By retaining the second-order terms, the SG system more accurately represents ageostrophic flows and frontal dynamics \cite{Cullen:2021}.

The SG system was introduced by Eliassen in 1949 \cite{Eliassen:1949} and later revisited by Hoskins in the 1970s \cite{Hoskins:1979}, who introduced \textit{geostrophic} coordinates. 
It has played a pivotal role in our understanding of large-scale (at length scales on the order of tens of kilometres) atmospheric dynamics and of atmospheric front formation.

The usefulness of the SG equations has been exemplified in operational settings, such as their use as a diagnostic tool by the UK Met Office, underscoring their value in meteorological \response{practice} \cite{Cullen:2018, Cullen:2021}.

Mathematically, this system came to prominence following the pioneering work of Brenier and Benamou \cite{Benamou1998WeakProblem}, who showed how the SG system in geostrophic variables is amenable to analysis using optimal transport techniques.
This formulation is the one we use in this paper as the basis of our numerical investigation. 

We use the SG system to model the evolution of an isolated cyclone, starting from a standard initial profile proposed by Sch\"ar and Wernli \cite{Schar:1993}.
Previous works, including those by Hoskins, Sch\"ar and Wernli, were constrained by the breakdown of the transformation between geostrophic and physical space in finite time, limiting the application of the SG system. 
In contrast, our approach leverages the optimal transport formulation, which overcomes this limitation and ensures that the transformation remains valid for all times, a fundamental advantage that extends beyond simply improving simulations. 
While existing studies focus on 2D computations of temperature and pressure evolution at the top and bottom boundaries, we go further by simulating the full 3D dynamics, capturing both the exterior surfaces and the interior of the domain.

The base state of this initial condition is a symmetric baroclinic jet combined with a uniform barotropic shear component controlled by the shear parameter, $A\in\R$, and given in terms of the non-dimensionalised pressure by
\begin{equation}\label{eq:SteadyState}
    \overline{\Phi}(\vb{x})=\frac{1}{2}\qty(\arctan(\frac{x_2}{1+x_3})-\arctan(\frac{x_2}{1-x_3}))-0.12x_2x_3-\frac{1}{2}A\qty(x_2^2-x_3^2).
\end{equation}
Importantly, this base state is harmonic and encodes the ramp-like structure observed in the presence of weather fronts.
The base state is then perturbed at the top and bottom of the domain via the perturbation function given by
\begin{equation}\label{eq:Perturbation}
    \begin{split}
    h(x_1,x_2)&= \qty(1+\qty(\frac{x_1}{0.5})^2+\qty(\frac{x_2}{0.5})^2)^{-\frac{3}{2}}-\frac{1}{2}\qty(1+\qty(\frac{x_1-1}{0.5})^2+\qty(\frac{x_2}{0.5})^2)^{-\frac{3}{2}} \\
    &\quad\quad-\frac{1}{2}\qty(1+\qty(\frac{x_1+1}{0.5})^2+\qty(\frac{x_2}{0.5})^2)^{-\frac{3}{2}}.
    \end{split}
\end{equation}
This perturbation is applied only to the temperature, i.e. to the derivative of the pressure with respect to the vertical coordinate, $x_3$. 

Our study seeks to advance the numerical treatment of the SG model by simulating its flow in 3D geostrophic space, continuing the research presented in \cite{Egan:2022} and iterated upon by \cite{Benamou:2024, Carlier:2024}. 
Specifically, we employ the damped Newton method recently developed by Kitagawa, Mérigot, and Thibert \cite{Kitagawa:2019} to evaluate numerically the semi-discrete transport map, which is equivalent to computing an optimal Laguerre tessellation of the 3D space. 
This method represented a significant advancement for numerical semi-discrete optimal transport methods and it aligns with our goal to adopt a mathematically rigorous and consistent formulation of the geometric method first proposed by Cullen and Purser \cite{Cullen:1984}. 
This approach is particularly desirable in view of its structural preservation qualities. 
Indeed, a crucial advantage of the semi-discrete optimal transport method over traditional finite element methods such as \cite{Visram:2014, Yamazaki:2017} is its capacity to preserve the underlying structures of the equations being discretised, so that solutions obtained through this method conserve total energy and simulate optimally mass-preserving flows within the fluid domain. 
Such characteristics are not only mathematically appealing but also crucial for the physical reliability of the simulations. 
This is particularly important when dealing with complex phenomena like frontal discontinuities, which are mathematically described as singularities occurring in finite time. 
Our numerical solutions, offer an accurate conservation of total energy, mirroring the physical behaviours observed in natural fluid dynamics and potentially allowing new insights into the understanding and prediction of atmospheric and oceanic phenomena.

\subsection{Semi-geostrophic system in discrete geostrophic variables}

In this section we provide the mathematical background for the model and explain, in brief, its connection to optimal transport. 
Consider a compact convex set $\X\subset \R^3$. $\X$ can be identified with the \textit{physical or fluid} domain. 
Furthermore, consider an open set $\Y\subseteq \R^3$, usually called \textit{geostrophic} space. 
The SG system in geostrophic space \cite{Benamou1998WeakProblem}, is given
\begin{equation}\label{mainPDE}
    \begin{split}
        &\partial_t \alpha_t + \mathrm{div}(\alpha_t \vb{v}[\alpha_t])=0, 
        \\
        &\vb{v}[\alpha_t]  = J(\mathrm{id}-T^{-1}),\qquad  
        J=\mqty(0 & -1 & 0 \\ 1 & 0 & 0 \\ 0 & 0 & 0),
    \end{split}
\end{equation}
where $\alpha:[0,\tau]\to\Pn(\Y)$ is a probability measure-valued map such that $\alpha_t=\alpha(t)\in\Pn(\Y)$ and $T:\X \to \Y$  is the optimal transport map from $\mathds{1}_{\X}$ (the normalised Lebesgue measure on $\X$) to  $\alpha_t$, the inverse of the potential vorticity \cite{Cullen:2021}, with respect to a quadratic cost. 
In this case the optimal transport map is defined as 
\begin{align*}
    T = \argmin_{\substack{T:\X \to \Y \\ T_{\#} \mathbbm{1}_{\X} = \alpha_t}} \int_{\X} \|\vb{x} - T(\vb{x})\|^2 \, \dd\vb{x},
\end{align*}
where $T_{\#} \mathbbm{1}_{\X}$ is the push-forward of the probability measure $\mathbbm{1}_{\X}$ under the map $T$ \cite{Santambrogio:2015}.
This is the case that has been studied most extensively, and for which there is a fairly exhaustive theory \cite{Santambrogio:2015}. 

Fix $N\in\mathbb{N}$. In our particle discretisation of this problem the target measure, $\alpha_t$, becomes the discrete probability measure $\alpha^N_t=\frac{1}{N} \sum_i^N \delta_{\vb{z}^i(t)}$, where $\vb{z}^i$ are points in $\Y\subset\R^3$. 
This yields an ODE which is given formally by
\begin{align}
    \dot{\vb{z}}^i &= \vb{v}[\alpha_t^N](\vb{z}^i), \nonumber
    \\
    \vb{v}\qty[\alpha_t^N](\vb{z}^i)  &=   J(\vb{z}^i-T^{-1}(\vb{z}^i)),
    \label{eq:FormalRHS} 
\end{align}
where $T:\X \to \Y$ is now the optimal transport map from $\mathbbm{1}_{X}$ to the discrete measure $\alpha_t^N$. 
This is a semi-discrete optimal transport problem, with respect to the quadratic cost. 
It is well known that its solution must be of the form $T=\frac{1}{N}\sum_{i=1}^N {\mathbbm 1}_{L^i}$, where $L^i$ are cells forming a covering of the space $\X$.
Rigorously, the $i$-th cell is defined, for $i\in\qty{1,\ldots,N}$, by
\begin{align*}
    L^i=\qty{\vb{x}\in\X : \|\vb{x}-\vb{z}^i\|^2-w^i \leq \|\vb{x}-\vb{z}^j\|^2- w^j \quad\forall\,j\in\qty{1,\ldots,N}},
\end{align*}
where $w^j\in\R$ guarantee that $\mathrm{vol}(L^i)=\frac{1}{N}$. 
Thus the solution of this transport problem is equivalent to solving for the {\em optimal tessellation} of the source space $\X$. 
For the quadratic cost, this is given by the so-called Laguerre tessellation \cite{Merigot:2021}. 
Note that Eq.~\eqref{eq:FormalRHS} is only formal, because $T$ is not invertible since its inverse could map single points to regions, $T^{-1}(\vb{z}^i)=L^i$.
Hence we replace $T^{-1}(\vb{z}^i)$ in Eq.~\eqref{eq:FormalRHS} by the centroid of the cell $L^i$. 
Then $\alpha^N$ is an exact weak solution of Eq.~\eqref{mainPDE} for the discrete initial data $\alpha^N_0$. 
The resulting ODE been studied extensively \cite{Bourne:2022} and is given by
\begin{equation}\label{mainODE}
    \vb{z}(t)=\qty(\vb{z}_1(t),\ldots,\vb{z}^N(t)), \qquad \dot{\vb{z}}=J^N\qty(\vb{z}-\vb{C}(\vb{z})),
\end{equation}
where $J^N\in\R^{3N\times 3N}$ is the block diagonal matrix $ J^N:=\mathrm{diag}(J,\ldots,J)$, and $\vb{C}(\vb{z})\in\X^{N}$ is the \emph{centroid map}, which identifies the centroid of each cell $L^i$ of the optimal tessellation:
\begin{align*}
    \vb{C}(\vb{z})=\qty(\vb{C}^1(\vb{z}), \ldots, \vb{C}^N(\vb{z})),\qquad C^i(\vb{z})=\frac 1 {\abs{L^i(\vb{z})}}\int_{L^i(\vb{z})}\vb{x} \,\dd\vb{x}.
\end{align*}
It is also energy conserving, where the total energy of the system is given by
\begin{align*}
    E(\vb{z}(t)) = \sum_{i=1}^N \int_{L^i}\|\vb{x}-\vb{z}^i(t)\|^2 \,\dd \vb{x}.
\end{align*}
Numerically this problem is solved efficiently, for a very large number of point particles, via the damped Newton method of \cite{Kitagawa:2019} already mentioned. 
We then select a 4th order Runge-Kutta (RK4) method as the one among existing ODE solvers for simulating the time evolution of the flow that achieves the best balance between efficiency and energy-conservation properties. 
This choice will be discussed further in \S~\ref{sec:results}.

\begin{rmk}
    It would be desirable to establish that the numerical solutions converge to a solution of the underlying partial differential equation (PDE) as the number of particles $N$ tends to infinity and the timestep size tends to zero. 
    The result in \cite{Bourne:2022} demonstrates convergence of a subsequence as $N\to\infty$, but this is not a full numerical analysis convergence result. 
    Moreover, convergence of the full sequence may be challenging to prove, as the uniqueness of weak solutions to the underlying PDE is not known. 
    In contrast, for the entropy-regularised fully discrete optimal transport scheme \cite{Carlier:2024}, convergence of the fully discrete scheme is established for $\varepsilon>0$, providing a stronger theoretical guarantee in this setting.
\end{rmk}

\subsection{Outline of the paper}

In Section \ref{sec:benchmark} we reproduce the 2D results of \cite{Egan:2022} and \cite{Benamou:2024} in order to validate and benchmark our code. 
In Section \ref{sec:initial} we present and discuss the 3D initial conditions used in our simulations which are informed by the work of \cite{Schar:1993} and are designed to generate an isolated large-scale tropical cyclone. 
We stress that, unlike \cite{Schar:1993} and other existing results that use these initial conditions to simulate the cyclone just on the surface of the domain, we compute numerically the evolution of the full 3D problem. 
In the last section we collect figures, showing both the evolution of the geostrophic particles and of the cells in physical space, that illustrate the evolution of the computed cyclonic flow. 
Finally in the appendix we show how the initial condition for the cyclonic flows is derived.

\section{2D Benchmark}\label{sec:benchmark}

\begin{figure}[!ht] 
    \includegraphics[scale=0.5]{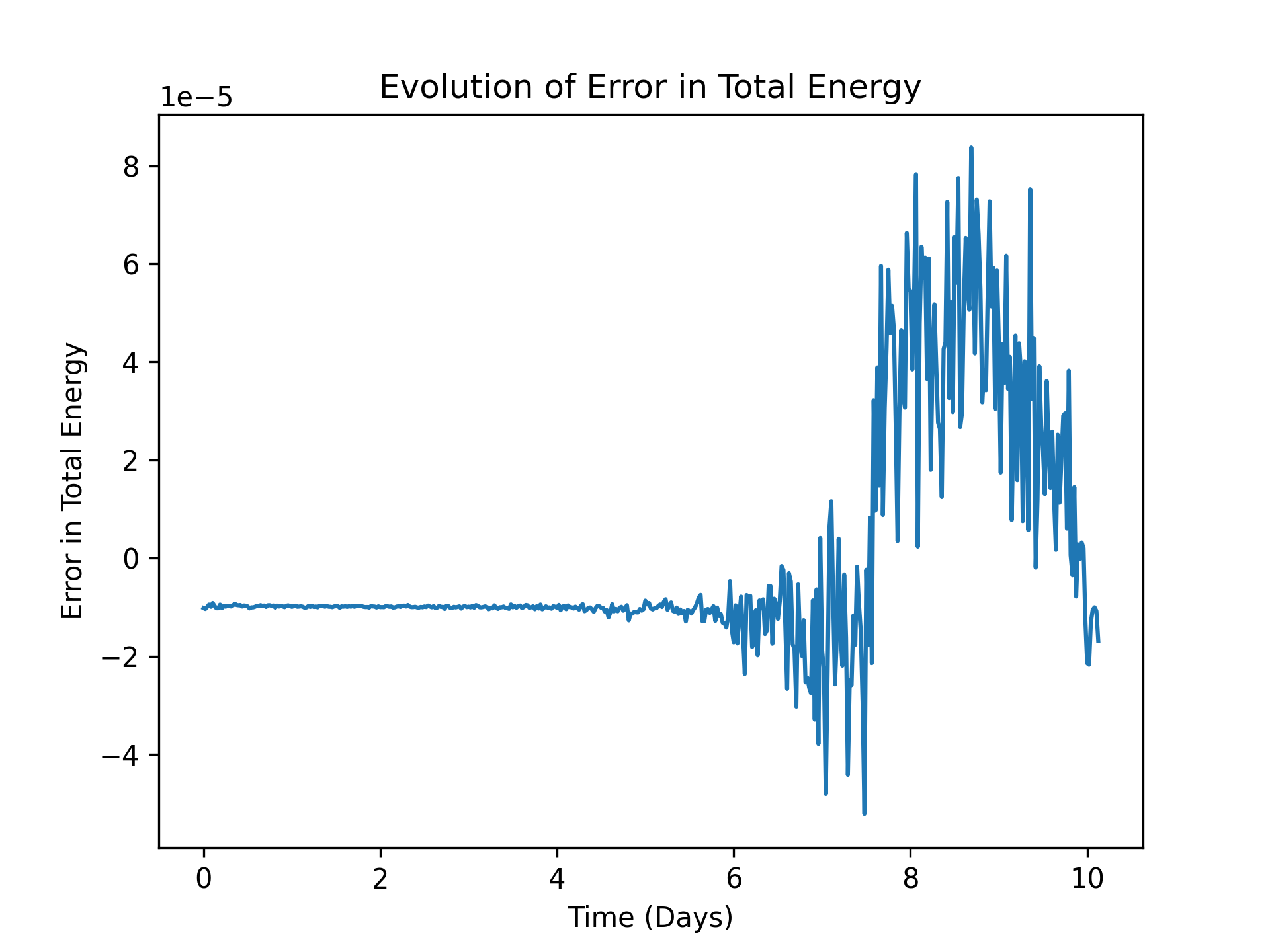}
    \caption{\response{The plot shows the evolution of the relative error in the total energy as defined in Eq.~\eqref{eq:errordefn}. The total energy fluctuates about $2.415\text{e}10$.}}
    \label{fig:2DConservation}
\end{figure}

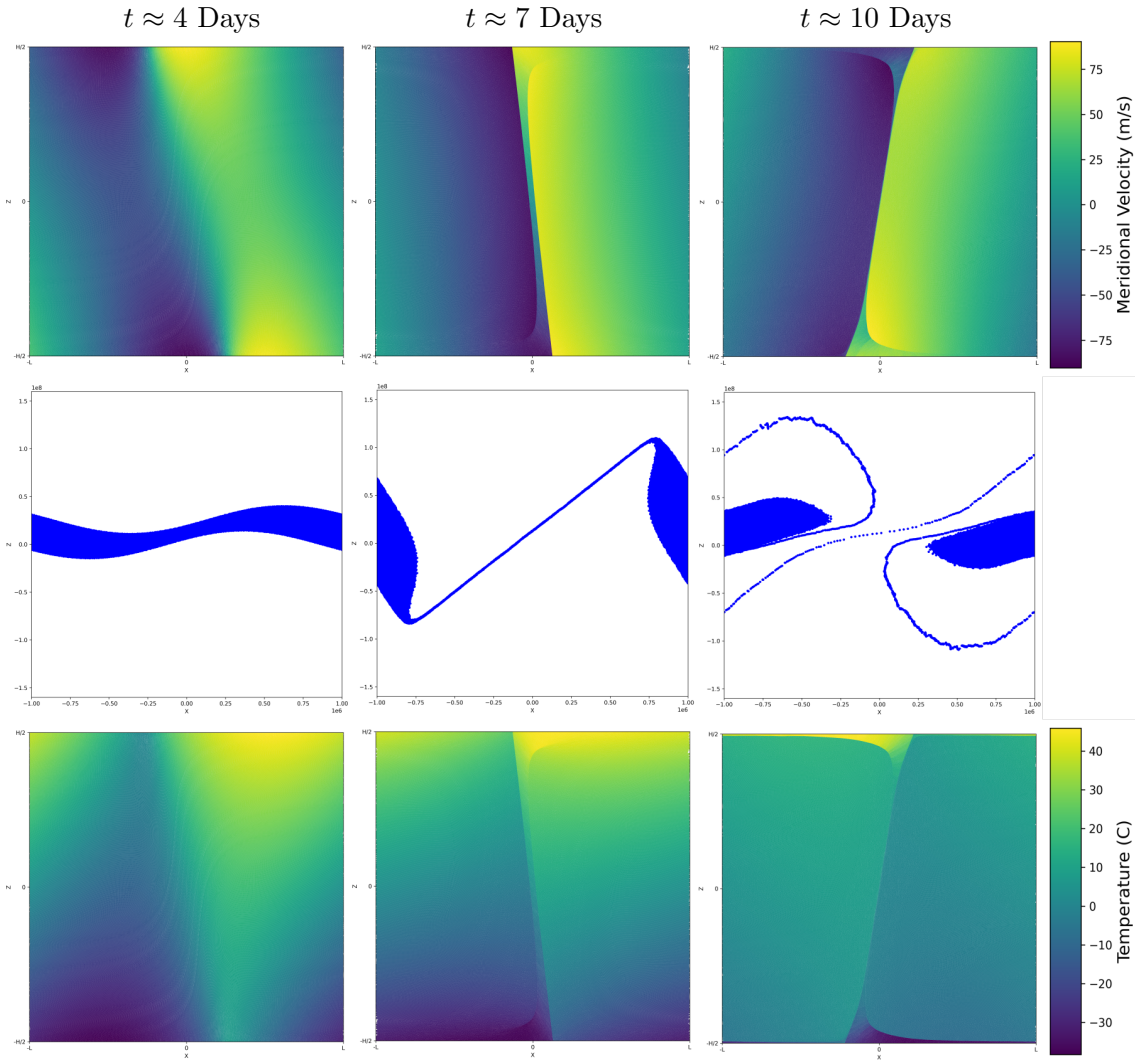
\begin{figure}[!ht]
    \centering
    \begin{tikzpicture}
        \node (img) at (0,0) {\includegraphics[scale=0.65]{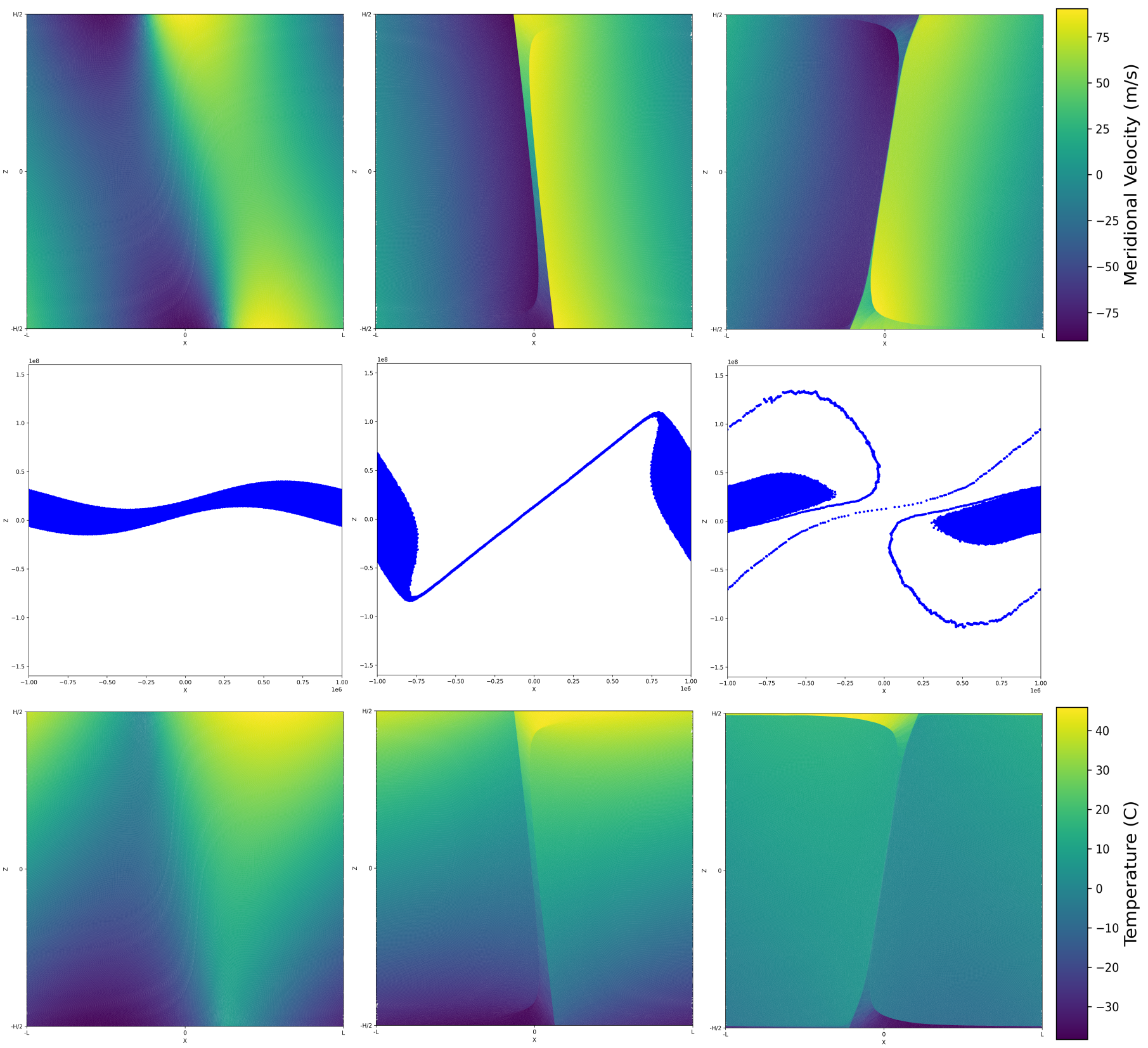}};
        \node at (-5,7) {$t\approx 4$ Days};
        \node at (-0.5,7) {$t\approx 7$ Days};
        \node at (4,7) {$t\approx 10$ Days};
    \end{tikzpicture}
    \caption{In the first row we display the evolution of the perturbation of the meridional velocity field in the physical space ($\X$). In the second row we display the evolution of the positions of the geostrophic particles (in $\Y$). Finally, in the third row we display the evolution of the perturbation of the temperature field in physical space. Over the course of 10 days we observe the formation and evolution of a weather front.}
    \label{fig:2DBenchmark}
\end{figure}

Egan et al.~in \cite{Egan:2022} first studied the application of a semi-discrete optimal transport scheme to solve the 2D incompressible SG equations. 
Given $\vb{z}(t)=(\vb{z}^1(t),\ldots,\vb{z}^N(t))\in\R^{2N}$ the ODE for the 2D system is 
\begin{equation*}
    \begin{cases}
        \dot{\vb{z}}^i= J(\vb{C}^i(\vb{z})-(\vb{z}^i\cdot\vu{e}_1)\vu{e}_1), \\
        \vb{z}^i(0)=\overline{\vb{z}}^i,
    \end{cases}
\end{equation*}
for $i\in\qty{1,\ldots, N}$, where $\overline{\vb{z}}^i$ denotes the initial position of the $i$-th seed, $\vb{c}^i$ is the centroid of the $i$-th 2D Laguerre cell, and 
\begin{align*}
    J=\mqty(0 & -1 \\ 1 & 0).
\end{align*}
The energy for the 2D system, up to a constant depending on the physical parameters of the problem, is given by
\begin{align*}
    E_{\textrm{2D}}(\vb{z}(t)) \simeq \frac{1}{2}\sum_{i=1}^N\qty( \int_{L^i} \|\vb{x}-\vb{z}^i(t)\|^2\,\dd\vb{x} - \qty(z_2^i(t)) ^ 2\int_{L^i}1\,\dd\vb{x}).
\end{align*}
Benamou et al.~\cite{Benamou:2024} built upon this foundation by implementing a fully discrete scheme, significantly enhancing simulation resolution. To validate our code, we replicated these results, employing Egan's technique in conjunction with Benamou's improved resolution. 
This enhancement was facilitated by advanced numerical schemes developed by Mérigot and Leclerc \cite{Leclerc2024Pysdot}.
We utilised the `unstable normal mode' scenario detailed in section 5.2 of \cite{Egan:2022} as a benchmark for our code. 
In contrast to \cite{Egan:2022}, whose highest resolution simulations were done with $N=2678$ particles, our simulations were done with $N=64284$ particles which is similar to the number of points used by \cite{Benamou:2024} ($N=65536$). 
As shown in Figure~\ref{fig:2DConservation}, our implementation conserves total energy with a relative error on the order of $10^{-5}$, comparable to the implementations in \cite{Egan:2022, Benamou:2024}. 
Egan et al.~used a timestep of 30 seconds and Benamou et al.~used a timestep of 91.44 minutes.
For our benchmark we employed a timestep of 30 minutes. 
Notably, with these advanced numerical schemes, we achieved a  relative error comparable to the one in \cite{Egan:2022} but with a timestep 60 times larger.
Furthermore, as illustrated in Figure~\ref{fig:2DBenchmark}, the system's evolution over the first 10 days aligns visually with the previous results of \cite{Egan:2022, Benamou:2024}, where we observe the formation of a weather front and its subsequent oscillations. 
In order to generate the plots we extract the meridional velocity ($v$) and temperature ($\theta$) from the seeds positions,
\begin{align*}
    v(\vb{x}, t) &= C_1\sum_{i=1}^N(\vb{z}_1^i(t)-\vb{C}_1^i(\vb{z},(t)))\mathbbm{1}_{L^i}\qty(\vb{x}) \\
    \theta(\vb{x}, t) &= C_2\sum_{i=1}^N\vb{z}_2^i(t)\mathbbm{1}_{L^i}\qty(\vb{x}),
\end{align*}
where $C_1, C_2 \in \R$ are physical constants (see \cite{Egan:2022}).

\newpage

\section{3D Benchmarks}\label{sec:initial}

To initialise our computations, we require a suitable initial condition.
Here we explain how we generate physical initial conditions for the 3D incompressible SG equations. 
The initial condition presented simulates the formation of an isolated large-scale tropical cyclone with or without initial shearing winds. 

\subsection{Initial condition generating an isolated semi-geostrophic cyclone}

In order to construct the initial condition for ODE (Eq.~\eqref{mainODE}) we need to find the initial seed positions $\overline{\vb{z}}$. 
These will be given by
\begin{align*}
    \overline{\vb{z}}^i = \qty(\mathrm{Id}+\widetilde{T})(\vb{x}^i),
\end{align*}
where the transport map, $T = \mathrm{Id} + \widetilde{T}$, is the gradient of a convex function, which is physically interpreted as the modified pressure. 
This means that the initial condition will be given by 
\begin{align*}
    \overline{\vb{z}}^i = \vb{x}^i + \nabla \Phi(\vb{x}^i),
\end{align*}
where the $\vb{x}^i$ are points in a uniform grid on the domain $\X$, and 
\begin{align*}
    \Phi(x_1,x_2,x_3)=\overline{\Phi}(x_1,x_2,x_3)+\widetilde{\Phi},(x_1,x_2,x_3)
\end{align*}
satisfies $\Delta\Phi=0$. 
Note that Laplace's equation arises from linearising of the Monge-Amp\`ere equation (see Appendix~\ref{sec:MA-LP} for more details). 
Thus in order to construct the initial condition we need to solve for the full pressure field. 
Notice, however, that the steady state modified pressure $\overline{\Phi}$, Eq.~\eqref{eq:SteadyState}, is harmonic. 
Therefore, to find the full pressure field we just need to find the effect that the perturbation of the temperatures on the surfaces has in the bulk of the domain. 
Thus we only consider $\widetilde\Phi$ that satisfies 
\begin{equation}\label{eq:ICSystem}
    \begin{cases}
        \Delta\widetilde{\Phi}=0, \\
        \widetilde{\Phi}(-a,\cdot,\cdot)=\widetilde{\Phi}(a,\cdot,\cdot), \\
        \widetilde{\Phi}(\cdot,-b,\cdot)=\widetilde{\Phi}(\cdot,b,\cdot), \\
        \pdv{\widetilde{\Phi}}{x_3}\eval_{x_3=0}=0.15h(x_1,x_2), \\\pdv{\widetilde{\Phi}}{x_3}\eval_{x_3=c}=-0.6h(x_1+1,x_2).
    \end{cases}
\end{equation}
The solutions of the system \eqref{eq:ICSystem} on a cuboid domain $\X=[-a,a]\times[-b,b]\times[0,c]$ approximates the effect that the surface perturbations have on the bulk. 
The solution satisfying the given boundary conditions, suitably adjusted for compatibility can be found explicitly (see Appendix~\ref{sec:PerturbSoln}). 
We use the values for $a$, $b$, and $c$ suggested in \cite{Schar:1993} : $c$ is the height of the lid set to be $0.45$ which corresponds to a physical height of 9 km, $a$ is 3.66 and $b$ is 1.75, which corresponds to a channel (periodic in $x$) of area $14640\times7000$ kilometers. 

\section{Results}\label{sec:results}

\subsection{Numerical Method}

The numerical method we employ consists of two main components, an optimal transport solver coupled with an ODE solver. 
For the optimal transport solver, we utilised the Pysdot package to generate Laguerre diagrams and solve the optimal transport problem using the damped Newton algorithm.
Detailed information on this approach can be found in \cite{Leclerc2024Pysdot, Merigot:2021, Egan:2022}. 
To enhance the stability and speed of convergence of the damped Newton algorithm, we applied a specific rescaling and translation of the initial configuration of geostrophic particles. 
Further details on this technique are provided in the work of \cite{Meyron:2019}. 
After obtaining the centroids of the Laguerre cells from the optimal transport solver, we applied a classical RK4 scheme to solve the ODE. 

\smallskip
\noindent The full code, including links to and instructions about dependencies, is publicly available and can be found at: \href{https://github.com/thelavier/3DIncompressibleSG}{\textbf{github.com/thelavier/3DIncompressibleSG}}

\subsubsection{Justification of ODE Solver}

\begin{figure}[!ht]
    \centering
    \begin{subfigure}[b]{0.48\textwidth}
        \includegraphics[width=\textwidth]{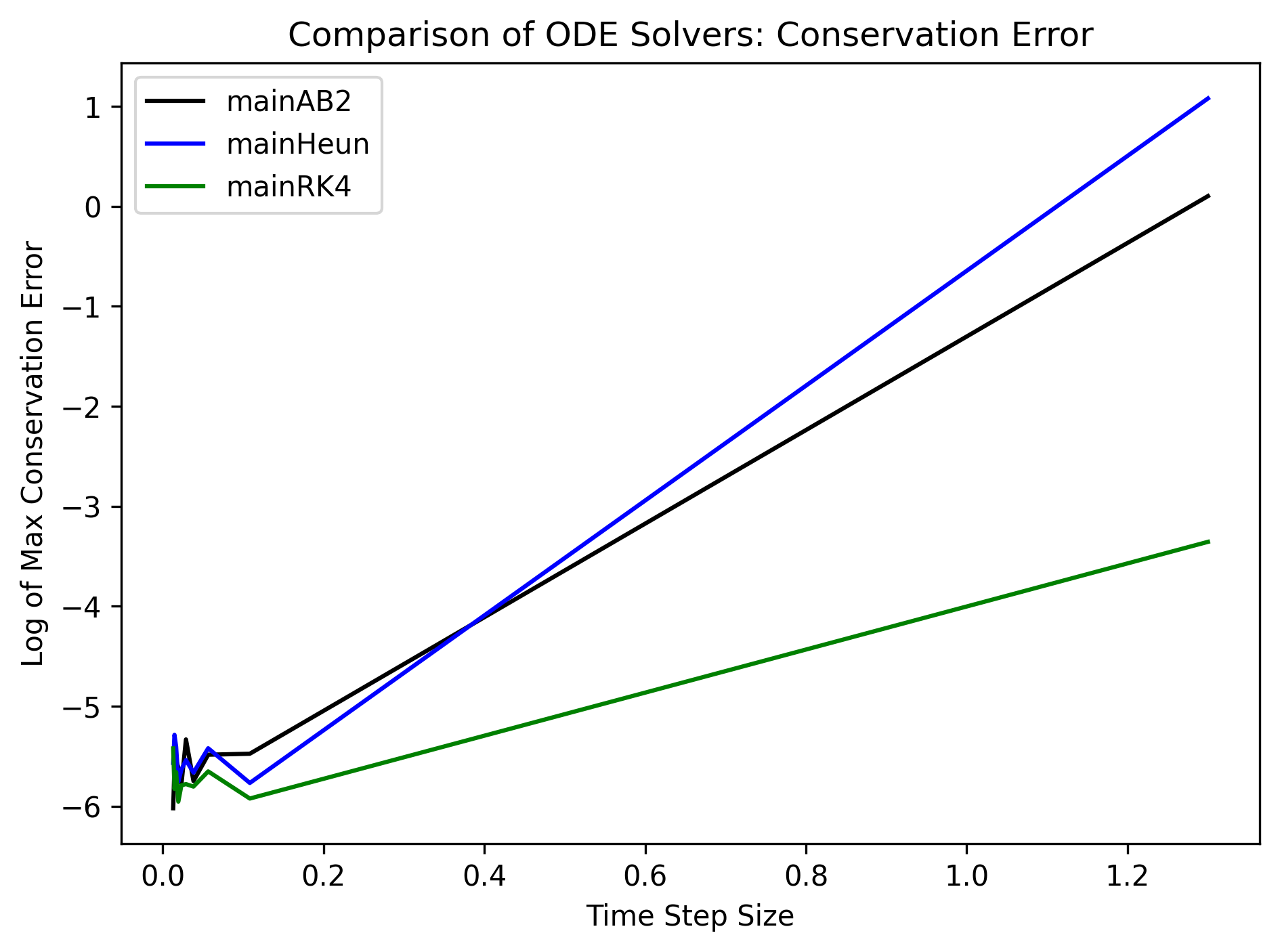}
        \caption{Here a comparison is presented showing how \response{the log of} the maximum relative conservation error in the energy changes with solver and step size.}
    \end{subfigure}
    \hfill 
    \begin{subfigure}[b]{0.48\textwidth}
        \includegraphics[width=\textwidth]{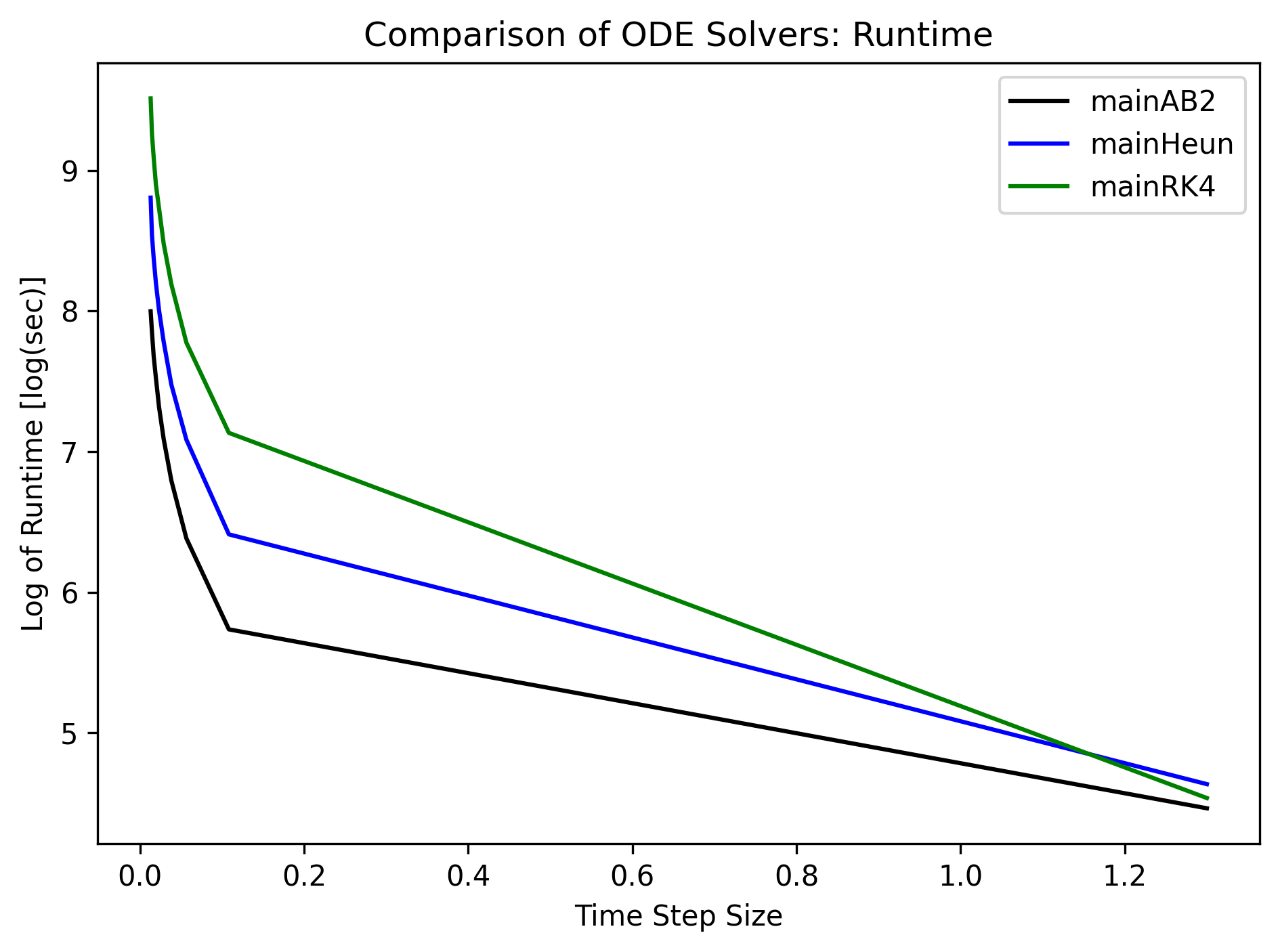}
        \caption{Here a comparison is presented showing how \response{the log of} the simulation runtime changes changes with solver and step size.}
    \end{subfigure}
    \caption{
    A comparison between AB2, Heun, and RK4 when coupled to an optimal transport solver. These two plots demonstrate the trade off between run time, time step size, and maximum relative error in the conservation of the energy. These plots support the idea that a balance can be struck between runtime, step size, and maximum relative error if one wants to run simulations in a ``reasonable" amount of time.}
    \label{fig:SolverComparisons}
\end{figure}

We experimented with several methods besides RK4, e.g. Adams-Bashforth 2 (AB2) and Heun. 
Following the lead of \cite{Benamou:2024} we ultimately chose RK4, even though it requires the solution of 4 optimal transport problems per timestep. 
We chose this method because of its energy conservation properties.
Ineed, while \cite{Benamou:2024} showed that there was no benefit in choosing a fourth order method over a second order method when considering the $W_1$ error with respect to a high resolution solution, we show that RK4 actually demonstrates a better performance over lower order methods in terms of energy conservation (see Figure \ref{fig:SolverComparisons}). 
This property enabled us to run 25-day simulations in 6-8 hours because, even though more optimal transport solves were required, the step size could be much larger without affecting the conservation of the energy, carrying forward the improvements observed in the 2D benchmark to the 3D simulations. 

\subsection{Experiments}\label{sec:experiments}

All the numerical experiments that we ran and their key parameters are presented in Table~\ref{tab:Error}.

\begin{table}[!ht]
\centering
\begin{tabular}{|cc|ccccccc|}
\hline 
& & \multicolumn{7}{|c|}{Wasserstein-2 Error} \\ 
& & $t\approx 2 \mathrm{~d}$ & $t\approx 4 \mathrm{~d}$ & $t\approx 9 \mathrm{~d}$ & $t\approx 13 \mathrm{~d}$ & $t\approx 17 \mathrm{~d}$ & $t\approx 21 \mathrm{~d}$ & $t\approx 25 \mathrm{~d}$ \\
\hline$\eta$ & 1 & 6.8667e-10 & 1.1666e-5 & 0.0068 & 0.0300 & 0.9804 & 0.4346 & 0.3678 \\
& 0.1 & 7.6761e-10 & 1.2229e-5 & 0.0064 & 0.0291 & 0.0871 & 0.2140 & 0.3373 \\
& 0.01 & 1.9691e-10 & 8.6135e-7 & 0.0038 & 0.0229 & 0.0740 & 0.1996 & 0.3042 \\
& 0.001 & - & - & - & - & - & - & - \\
\hline$h[\mathrm{sec}]$ & 10803.58 & 0.0091 & 0.0223 & 0.0866 & 0.3650 & 0.6162 & 0.5932 & 0.9808 \\
& 7190.46 & 0.0022 & 0.0071 & 0.0393 & 0.1865 & 0.5765 & 0.7715 & 0.8195 \\
& 5388.39 & 0.0014 & 0.0055 & 0.0428 & 0.1645 & 0.6625 & 0.9556 & 0.8565 \\
& 3595.23 & 0.0007 & 0.0032 & 0.0351 & 0.1377 & 0.7114 & 0.9401 & 1.0594 \\
& 2700.90 & 0.0004 & 0.0025 & 0.0331 & 0.1433 & 0.3229 & 0.5205 & 0.4602 \\
& 1799.11 & 5.2552e-5 & 0.0017 & 0.0246 & 0.0939 & 0.4740 & 0.9702 & 0.6124 \\
& 899.93 & - & - & - & - & - & - & - \\
\hline$N$ & 4096 & 0.1600 & 0.2594 & 0.6186 & 1.3901 & 1.8708 & 1.9466 & 1.4330 \\
& 5832 & 0.1233 & 0.1503 & 0.4656 & 1.0459 & 1.4215 & 1.3115 & 0.8541 \\
& 10648 & 0.0699 & 0.0939 & 0.3040 & 0.5674 & 0.9495 & 1.2798 & 0.9993 \\
& 21952 & 0.0246 & 0.0315 & 0.1190 & 0.7100 & 0.7642 & 1.4379 & 0.8179 \\
& 32768 & - & - & - & - & - & - & - \\
\hline
\end{tabular}
\caption{\response{Wasserstein-2 error, as defined by Eq.~\eqref{eq:WassError}, between simulations at the highest resolution (reference solution) and lower-resolution simulations as parameters vary. This table specifically presents errors resulting from the ODE solver used for temporal integration, rather than the damped Newton solver used for the optimal transport problem. The error dependence on the solver tolerance ($\eta$), timestep size ($h$), and particle count ($N$) is shown.} Simulations investigating the impact of $\eta$ were done with $h=3595.23$ and $N=32768$. Simulations investigating the impact of $h$ were done with $\eta=0.001$ and $N=32768$. Simulations investigating the impact of $N$ were done with $\eta=0.001$ and $h=3595.23$.}
\label{tab:Error}
\end{table}

\subsubsection{Error Computation}

In order to analyse the quantitative performance of our method, we consider the following error
\begin{equation}\label{eq:WassError}
    \mathrm{Error}(t) = W_2^2\qty(\frac{1}{N}\sum_{i=1}^N\delta_{\vb{z}^i_{\mathrm{true}}(t)}, \frac{1}{N}\sum_{i=1}^N\delta_{\vb{z}^i_{\mathrm{approx}}(t)}).
\end{equation}
\response{
Here $\vb{z}_{\mathrm{true}}(t)$ is the ensemble of seed positions generated by our finest-resolution run, which we adopt as a surrogate ``ground truth," while $\vb{z}^i_{\mathrm{approx}}(t)$ are the corresponding seeds produced by the lower-resolution simulation whose accuracy we wish to evaluate. 
The squared 2-Wasserstein distance between the two empirical measures therefore quantifies, in a physically meaningful way, how far the approximate Lagrangian configuration deviates from the reference at time $t$.
}
This allows us to analyse the effect of changing the three key simulation parameters : percent tolerance of the optimal transport solver, time step size, and number of particles. 
We are also interested the the ability of the solver to preserve averaged quantities, in particular the energy. 
In order to investigate this we also considered
\begin{equation}\label{eq:errordefn}
    \text{Error in Total Energy}(t) = \frac{E_{\mathrm{mean}}(\vb{z}) - E(\vb{z}(t))}{ E_{\mathrm{mean}}(\vb{z}) },
\end{equation}
and
\begin{align*}
    \text{Max Conservation Error} = \max_{t\in[0,\tau]} \text{Error in Total Energy}(t).
\end{align*}
In order to compute efficiently the Wasserstein-2 error for our analysis, we employed Jean Fedey's GeomLoss package to compute the Sinkhorn divergence approximation of the Wasserstein-2 distance \cite{feydy2019interpolating}. 
The Sinkhorn divergence approximation of the Wasserstein-2 distance is an efficient way of measuring the distance between two probability measures.
All errors were computed against a ``ground-truth" simulation, defined as the highest resolution simulation. 
All simulations were conducted using the RK4 method. 

\subsubsection{Results of Experiments}

\begin{figure}[!ht]
    \centering
    \includegraphics[scale=0.5]{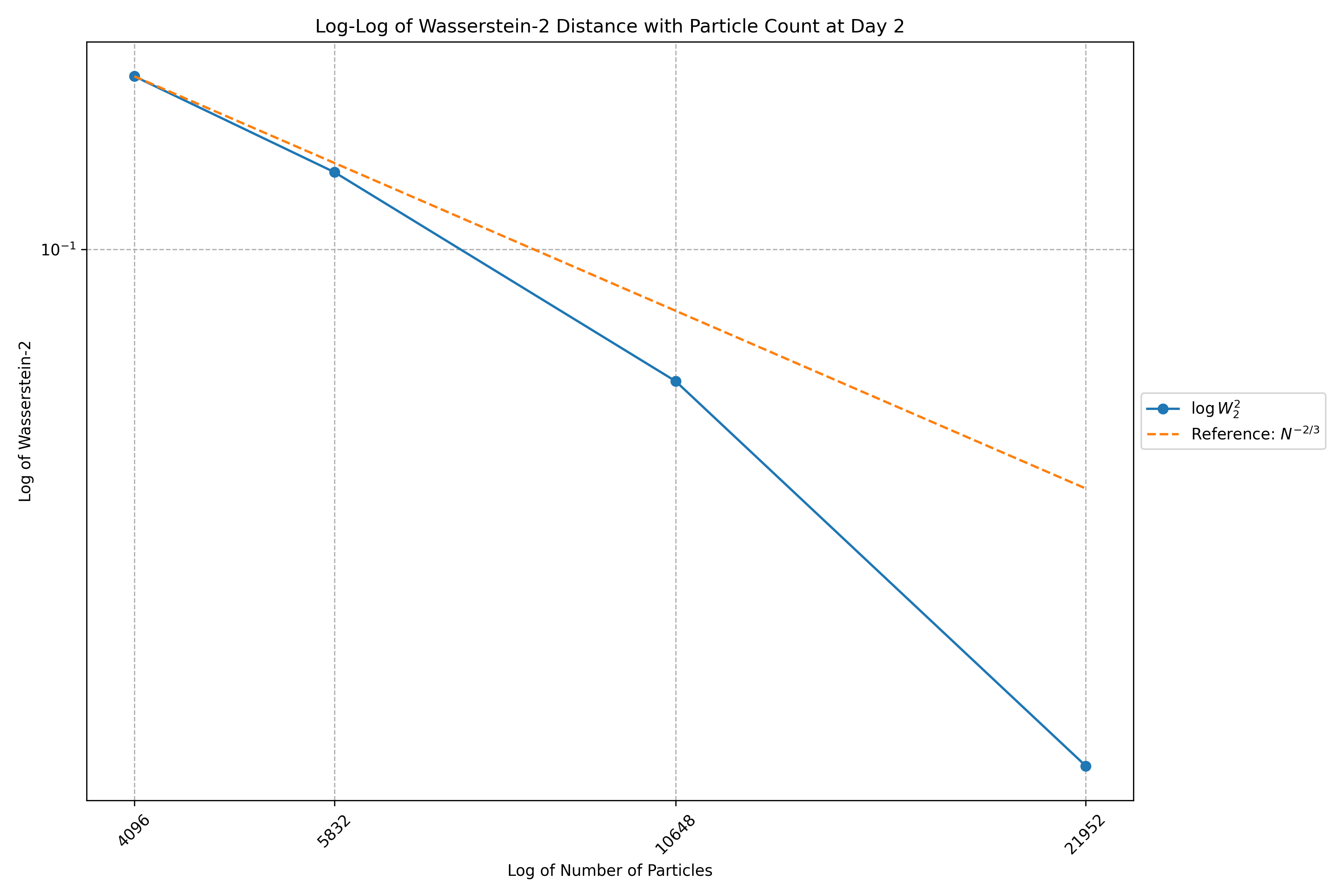}
    \caption{Log-log plot of the change in the Wasserstein-2 error at day 2 with respect to the change in the number of particles (in blue). In orange is a plot of the theoretical best decrease in the discretisation error with respect to the number of particles.}
    \label{fig:Particles}
\end{figure}

\begin{figure}[!ht]
    \centering
    \includegraphics[scale=0.5]{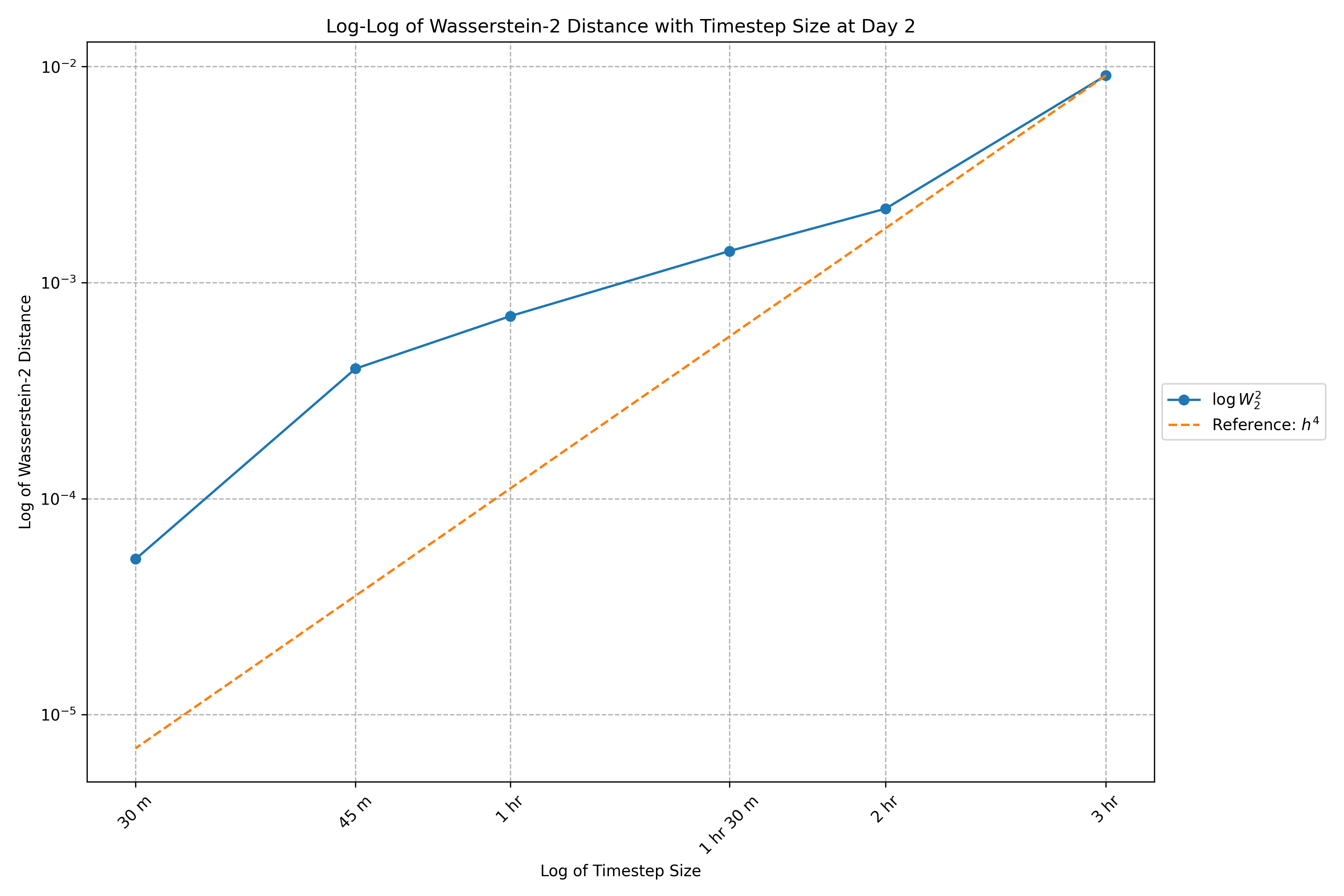}
    \caption{Log-log plot of the change in the Wasserstein-2 error at day 2 with respect to the change in the size of the timestep (in blue). In orange is a plot of the theoretical best decrease in the error with respect to the timestep size for Runga-Kutta 4.}
    \label{fig:Timestep}
\end{figure}

\begin{figure}[!ht]
    \includegraphics[scale=0.5]{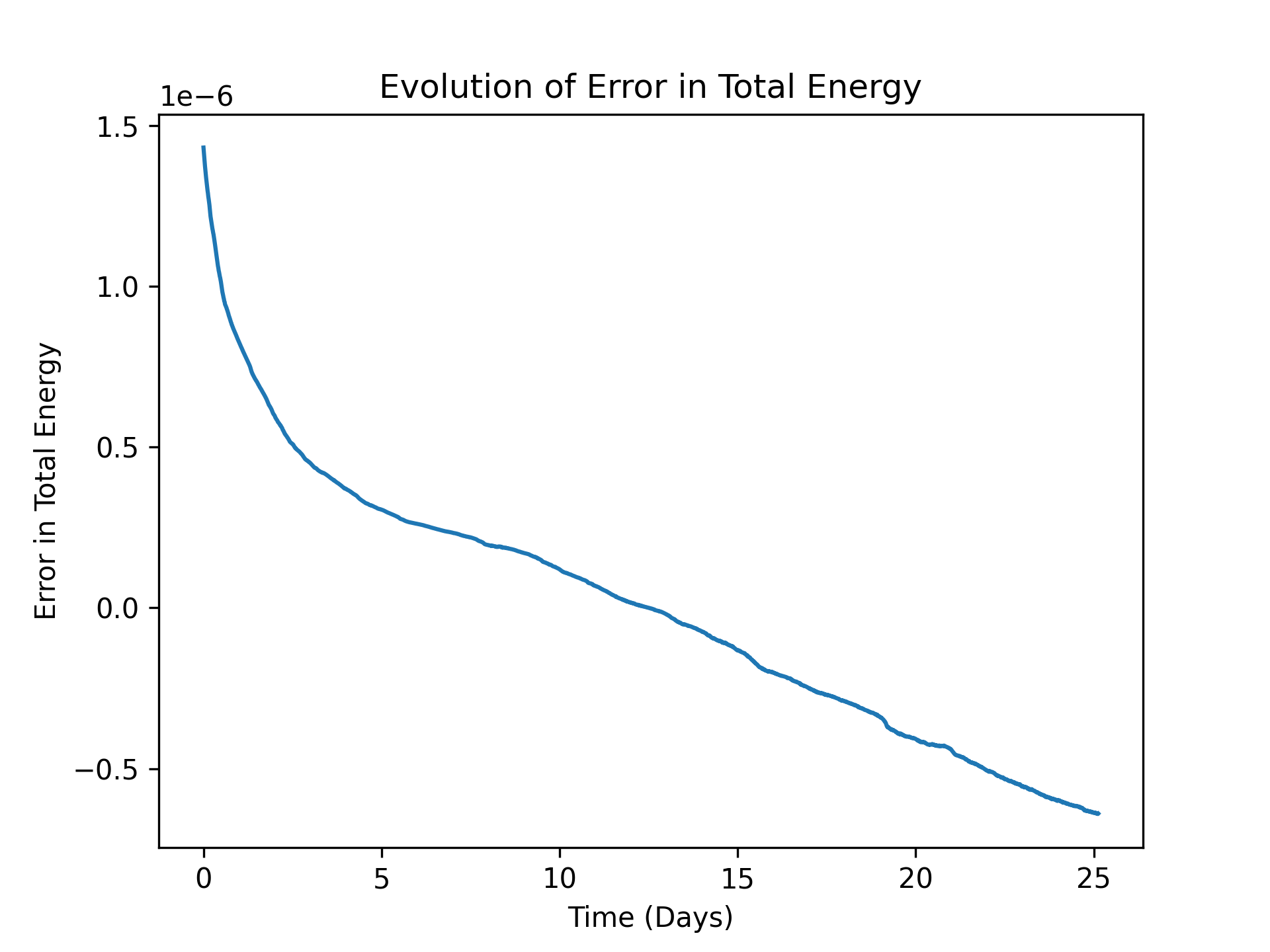}
    \caption{\response{The plot shows the evolution of the relative error in the total energy as defined in Eq.~\eqref{eq:errordefn}. The total energy fluctuates about $3.799$.}}
    \label{fig:3DConservation}
\end{figure}

As shown in Table~\ref{tab:Error}, for short durations such as day 2 and 4, we observe the expected reduction in error with respect to timestep size $\qty(h^4)$ and particle count $\qty(N^{-\flatfrac{2}{3}})$ \cite{Kloeckner:2012}, as illustrated in Figures~\ref{fig:Particles} and \ref{fig:Timestep}. 
However, as the simulation extends to 25 days, we observe a deterioration in the accuracy, and the anticipated decay in error relative to the number of particles and timestep size no longer holds. 
This decline is not unexpected, as the solution may become singular in physical space, where standard truncation error estimates only apply to smooth solutions. 
The lack of significant change in the error when varying the tolerance on the optimal transport solver is also unsurprising, as the damped Newton method often overshoots the specified tolerance. 
Consequently, within reasonable bounds, the choice of tolerance for the optimal transport solver has a limited impact on the overall accuracy. 
\response{Averaged properties of the system, such as total energy, remain well conserved, as shown in Figure~\ref{fig:3DConservation}, indicating that while the accuracy degrades, the overall physical integrity of the simulation is preserved.}
To address issues with the accuracy of the \response{ODE} solver, we plan to follow up with a study of the spectral properties of the system, particularly focusing on the centroid map, to better understand this deterioration and develop methods that maintain higher accuracy over long-term simulations. 
Preliminary analysis suggests that the issue may be related to a large relative spectral gap in the eigenvalues of the centroid map.

\newpage

\subsection{Observations of Cyclones}

In what follows is a qualitative discussion of the evolution of an isolated large-scale tropical cyclone as controlled by the SG system.
To aid in this discussion we extract the zonal velocity (ZVel), meridional velocity (MVel), \response{geostrophic velocity}, total geostrophic velocity (TVel), and temperature :
\begin{align*}
    \mathrm{ZVel}(\vb{x}, t) &= \sum_{i=1}^N \qty(C_2^i(\vb{z}(t)) - z_2^i(t)) \mathbbm{1}_{L^i}\qty(\vb{x}), \\
    \mathrm{MVel}(\vb{x}, t) &= \sum_{i=1}^N \qty(z_1^i(t) - C_1^i(\vb{z}(t))) \mathbbm{1}_{L^i}\qty(\vb{x}), \\
    \response{\vb{u}_g(\vb{x},t)} &\response{= \qty(\mathrm{ZVel}(\vb{x}, t), \mathrm{MVel}(\vb{x}, t))^T} \\
    \mathrm{TVel}(\vb{x}, t) &= \sum_{i=1}^N \norm{\vb{u}_g(\vb{x},t)} \mathbbm{1}_{L^i}\qty(\vb{x}), \\
    \text{Temperature}(\vb{x}, t) &= \sum_{i=1}^N z_3^i(t) \mathbbm{1}_{L^i}\qty(\vb{x}).
\end{align*}
We also compute the root mean squared velocity (RMSv) of the three different velocities 
\begin{align*}
    \mathrm{RMSv} = \sqrt{\frac{1}{\X}\int_{\X}\abs{v(\vb{x},t)}^2 \,\dd \vb{x}},
\end{align*}
to support our consideration of different initial shearing regimes.

\begin{figure}[!ht]
    \centering
    \begin{tikzpicture}
        \node (img) at (0,0) {\includegraphics[scale=0.375]{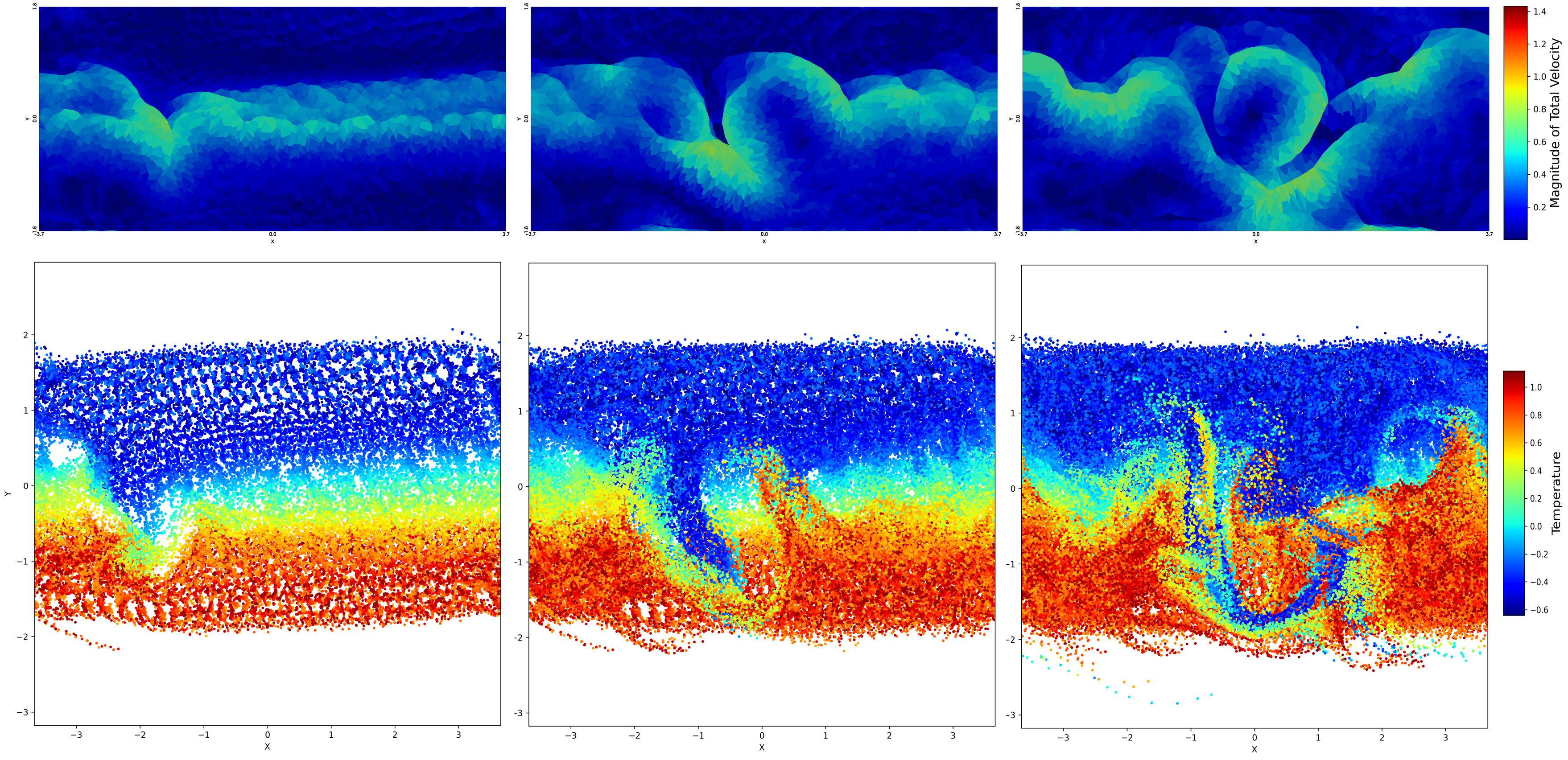}};
        \node at (-5.5,4.5) {$t\approx 4$ Days};
        \node at (-0.25,4.5) {$t\approx 8$ Days};
        \node at (5,4.5) {$t\approx 12$ Days};
    \end{tikzpicture}
    \hspace{2cm}
    \begin{tikzpicture}
        \node (img) at (0,0) {\includegraphics[scale=0.375]{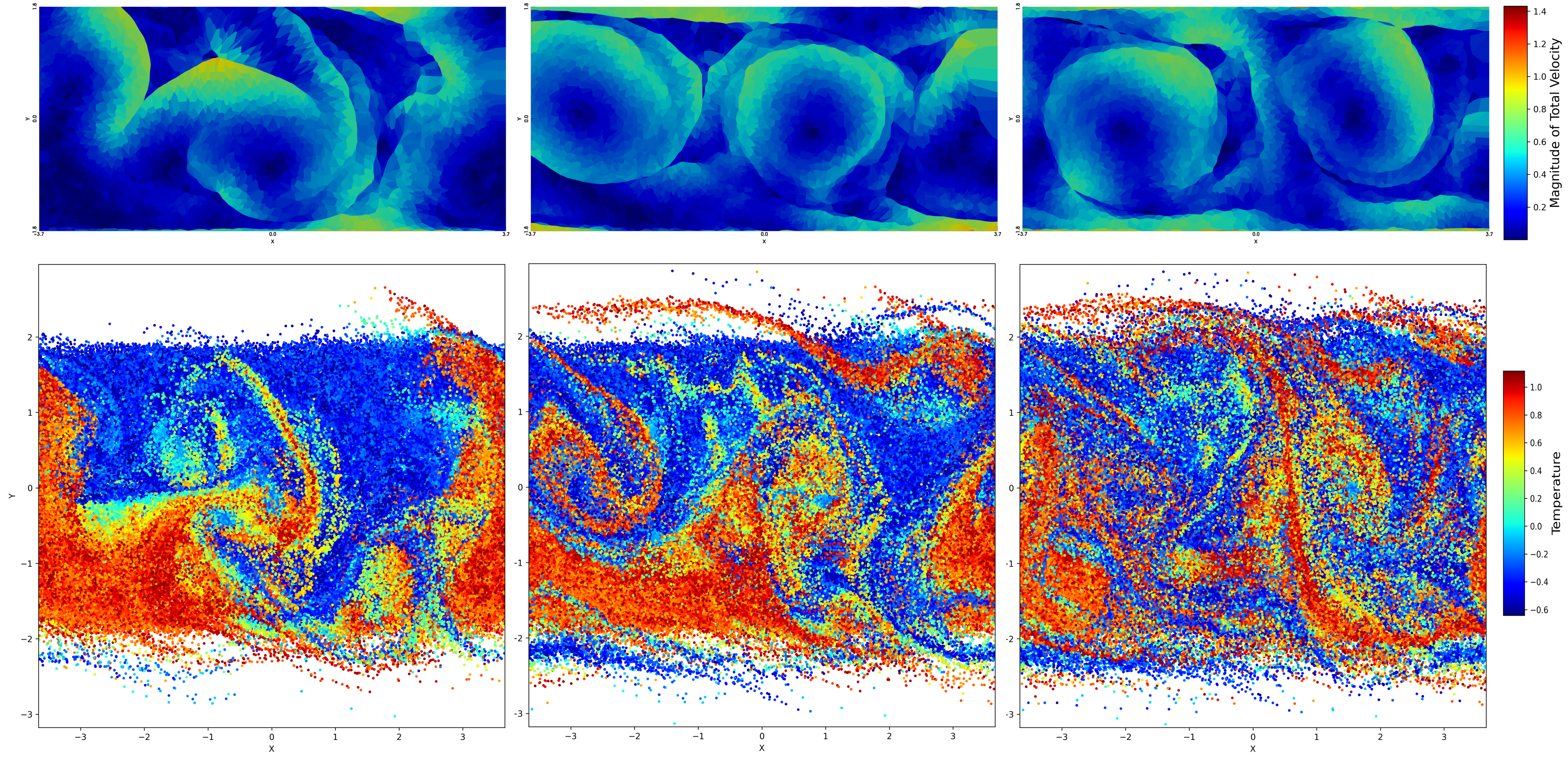}};
        \node at (-5.5,4.5) {$t\approx 16$ Days};
        \node at (-0.25,4.5) {$t\approx 20$ Days};
        \node at (5,4.5) {$t\approx 25$ Days};
    \end{tikzpicture}
    \caption{All images are done with the camera looking down on the top of the domain. In rows one and three we see the evolution of the magnitude of the total geostrophic velocity (in physical space, $\X$) over 25 days. In rows 2 and 4 we see the evolution of the position of the geostrophic particles over 25 days. In geostrophic space temperature corresponds to position in the third dimension. In these images this is captured in the colouring of the particles. Colder  ``lower" particles are blue and hotter ``higher" particles are red. Simulation parameters: $N=64000$, $\eta=10^{-3}$, and $h=30$ min.}
    \label{fig:CycloneEvolution}
\end{figure}
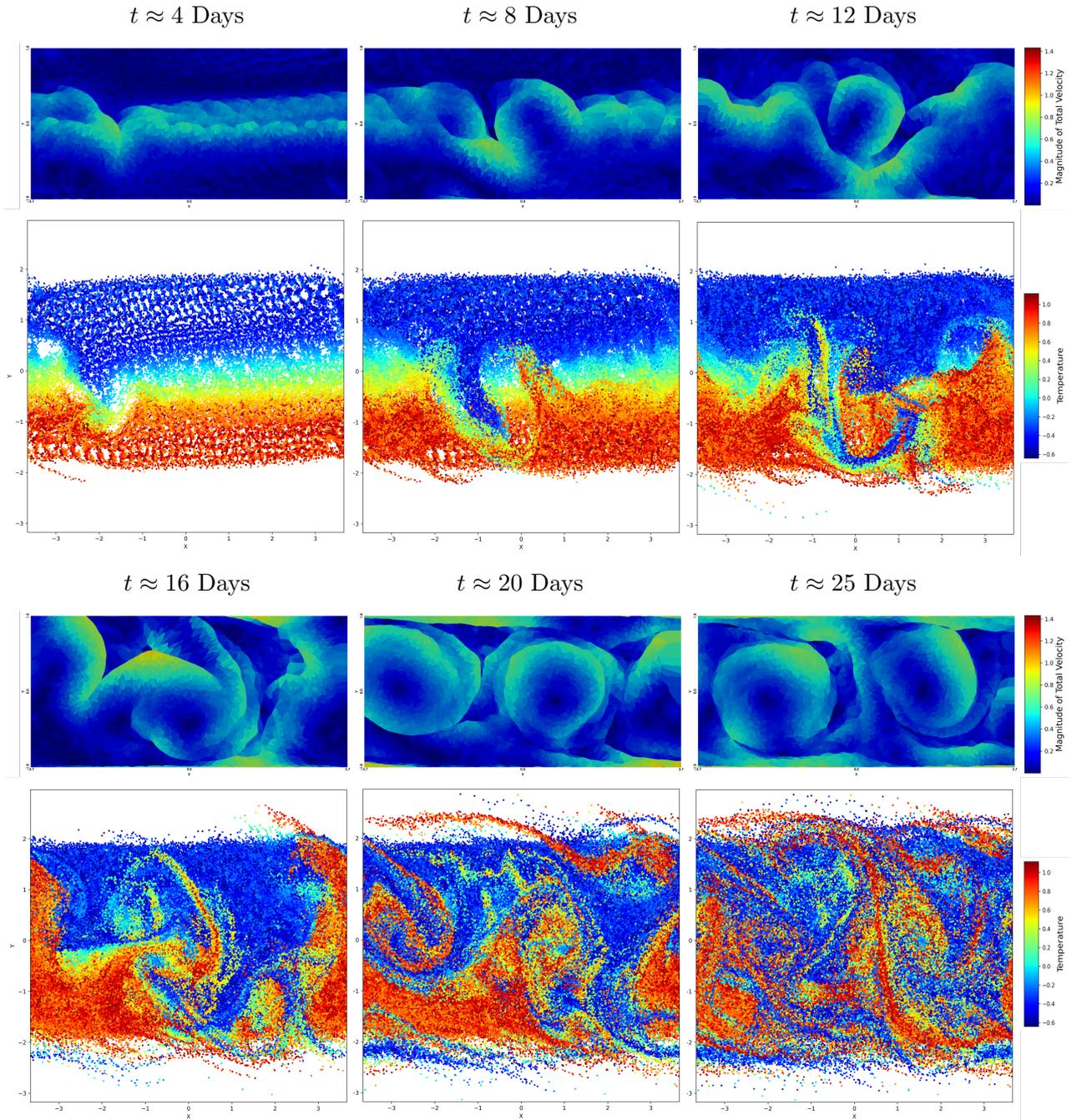

In Figure~\ref{fig:CycloneEvolution}, we observe the evolution of a 3D incompressible SG system where cold and hot air masses are initially separated and subsequently mix.
The series of images track the development of this interaction over a period of 25 days. 
In rows one and three, the images display the magnitude of the total geostrophic velocity in physical space ($\X$). 
Initially, at $t\approx4$ days, a distinct front forms between the cold and hot air masses. 
As time progresses to $t\approx8$ and $t\approx12$ days, an instability along this front propagates, evolving into a chain of rotational systems, indicative of cyclone and anticyclone formation.

Rows two and four depict the evolution of seed positions in geostrophic space ($\Y$), where the temperature corresponds to the vertical position in the third dimension. 
The geostrophic particles are color-coded to represent their vertical positions: blue for colder or ``lower" and red for hotter or ``higher". 
At $t\approx4$ days, the seeds are relatively evenly distributed along the front. 
By $t\approx8$ and $t\approx12$ days, the seeds begin to cluster and spiral, showing the development of vortices as the system becomes more dynamic.

By $t\approx16$ days, the rotational structures become more pronounced, and by $t\approx20$ and $t\approx25$ days, the system displays fully developed cyclonic and anticyclonic patterns. 
This visual evidence supports the conclusion that the initial instability evolves into a series of complex, rotating systems. 
The continued development and interaction of these vortices demonstrate the non-linear and chaotic nature of the 3D incompressible SG dynamics, while energy is conserved throughout the process, as indicated by the stable total energy observed in long-term simulations.

\begin{figure}[!ht]
    \centering
    \begin{tikzpicture}
        \node (img) at (0,0) {\includegraphics[scale=0.5]{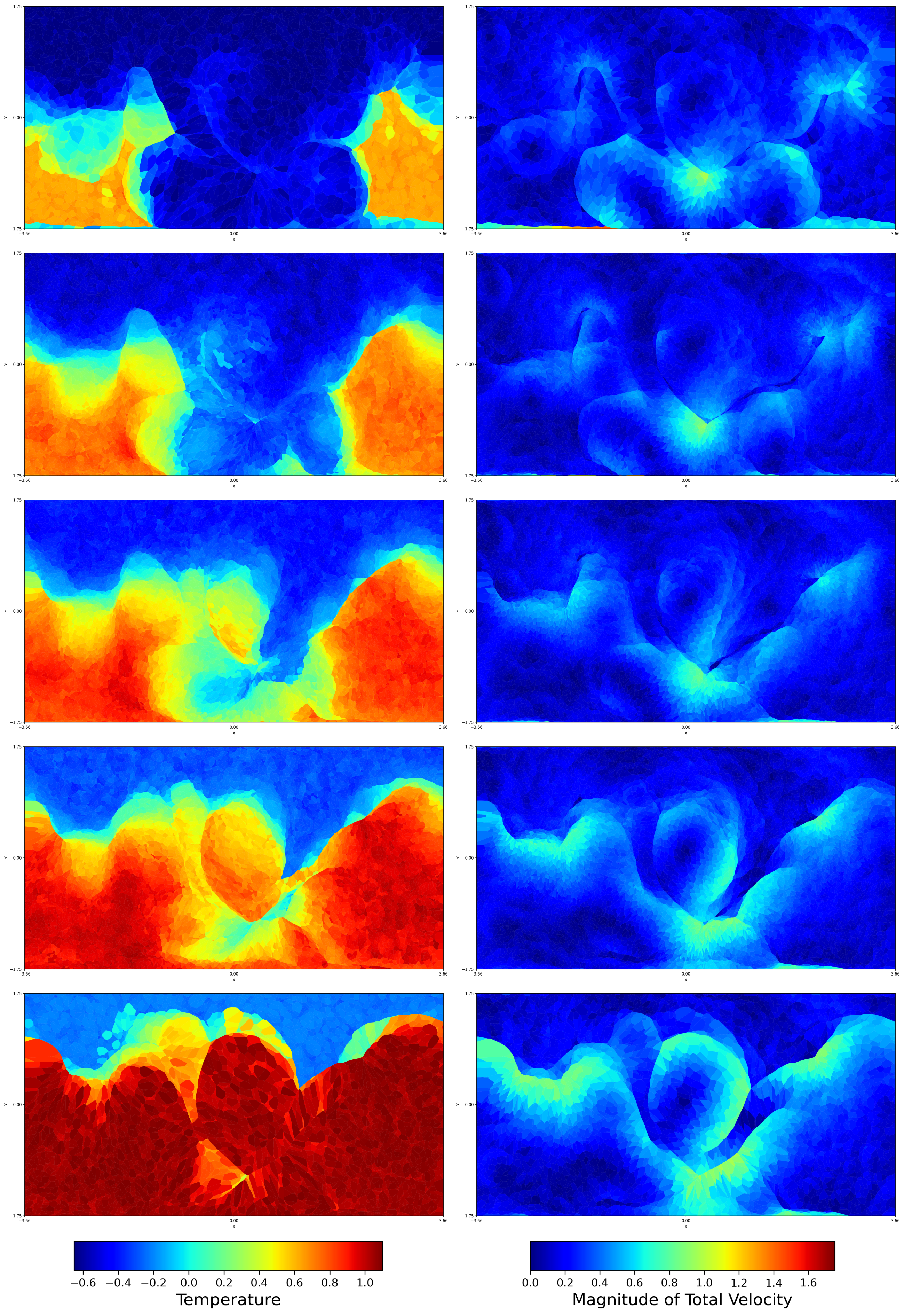}};

        \node[rotate=90] at (-5.75, 6.5) {\small$h = 0 \,$ km};
        \node[rotate=90] at (-5.75, 3.5) {\small$h = 2.25\,$ km};
        \node[rotate=90] at (-5.75, 0.5) {\small$h = 4.5\,$ km};
        \node[rotate=90] at (-5.75, -2.5) {\small$h = 6.75\,$ km};
        \node[rotate=90] at (-5.75, -5.5) {\small$h = 9\,$ km};
    \end{tikzpicture}
    \caption{Vertical slices of temperature (left column) and total velocity magnitude (right column) at various altitudes ($h = 0$ km, 2.25 km, 4.5 km, 6.75 km, and 9 km) after 12 days of simulation. The temperature distribution highlights regions of significant thermal activity, with warmer areas denoted by red hues and cooler areas by blue, indicating the presence of convective structures and stratification. The velocity magnitude plots reveal the structure of the flow, with areas of higher velocity depicted in yellow-green, illustrating the dynamics of the developing cyclone and associated turbulence. These slices offer a detailed view of the interaction between thermal and kinematic fields throughout the bulk of the cyclone, emphasising the formation and behaviour of flow structures across multiple atmospheric layers. Simulation parameters: $N=64000$, $\eta=10^{-3}$, and $h=30$ min.}
    \label{fig:Slices}
\end{figure}
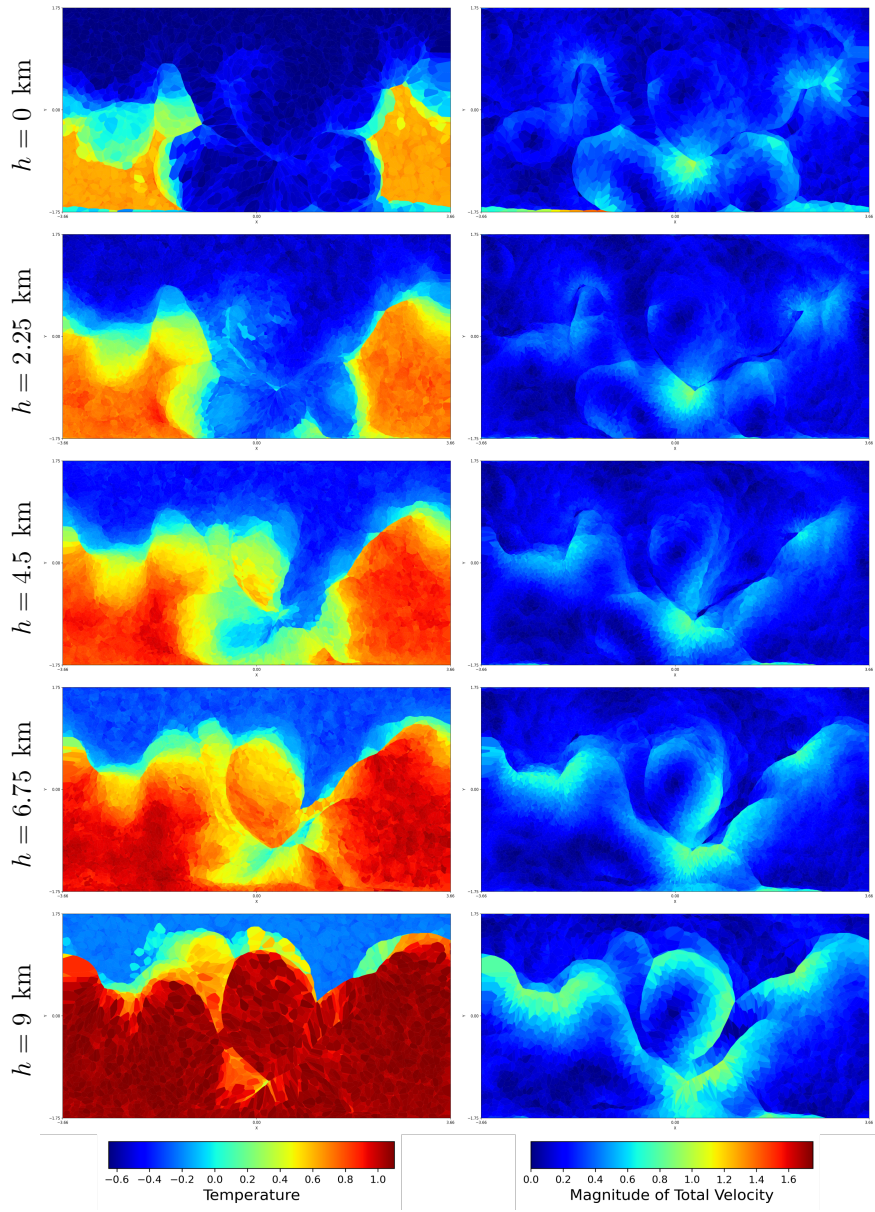

Figure~\ref{fig:Slices} complements this analysis by presenting horizontal cross-sections of the temperature and velocity magnitude at different altitudes within the domain after 12 days of evolution. 
These slices show the vertical structure of the flow, revealing the coupling between thermal and dynamical processes. 
The interaction between cold and warm air masses drives the development of baroclinic instability, leading to the characteristic rotation observed in geostrophic systems.

\newpage

\subsubsection{Shear Parameter}

Finally, as demonstrated in the seminal works by Davies et al.~\cite{Davies1991TheDevelopment} and Wernli et al.~\cite{Wernli1998TheSimulations}, the impact of a shearing wind on cyclone development is significant, with the evolution of the system being highly sensitive to the horizontal shear imposed at the initial time. 
These previous studies explored how variations in background shear influence the formation and characteristics of cyclonic and frontal structures. 
A key limitation they faced, however, was the loss of regularity of the transformation between geostrophic and physical space occurring between 4 to 8 days, restricting their ability to investigate the long-term effects of the shear.

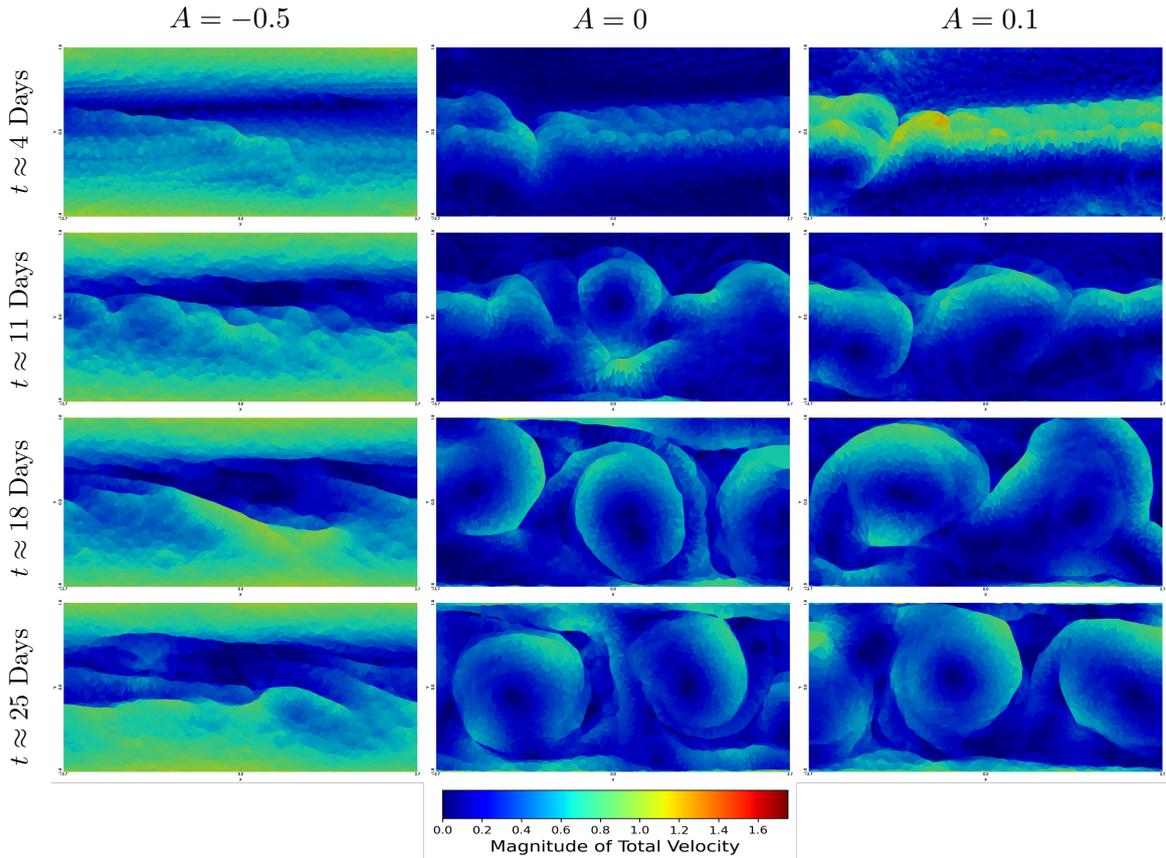
\begin{figure}[!ht]
    \centering
    \begin{tikzpicture}
        \node (img) at (0,0) {\includegraphics[scale=0.35]{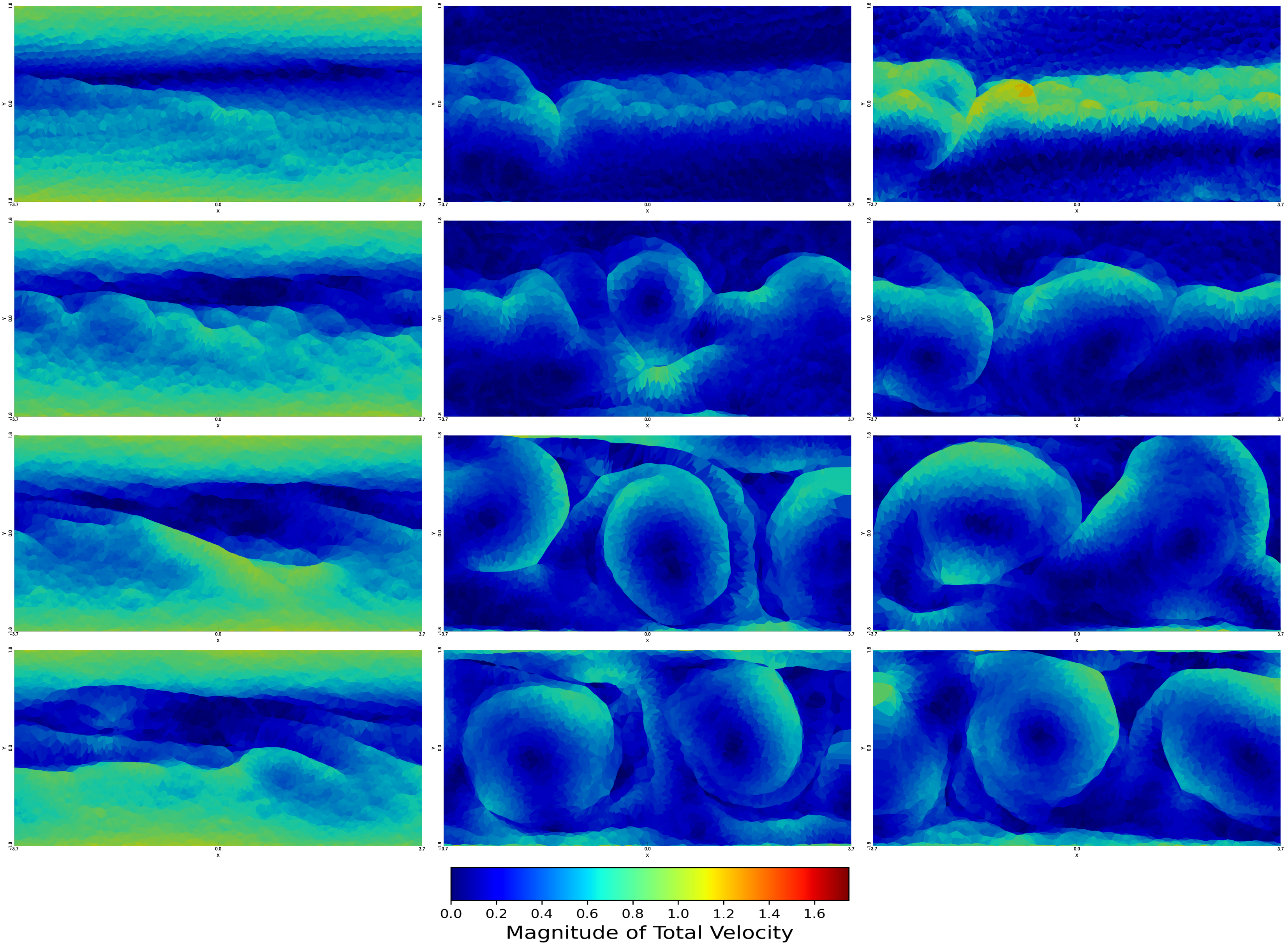}};
        \node at (-5,5.75) {$A=-0.5$};
        \node at (0,5.75) {$A=0$};
        \node at (5,5.75) {$A=0.1$};

        \node[rotate=90] at (-7.75, 4.25) {\small$t \approx 4$ Days};
        \node[rotate=90] at (-7.75, 1.75) {\small$t \approx 11$ Days};
        \node[rotate=90] at (-7.75, -0.75) {\small$t \approx 18$ Days};
        \node[rotate=90] at (-7.75, -3.25) {\small$t \approx 25$ Days};
    \end{tikzpicture}
    \caption{Here we demonstrate how adding a background shearing wind to the initial condition results in the disruption of or strengthening of cyclone formation over 25 days. In the first column we see the effect of a strong shear wind completely disrupting cyclone formation in the channel, in contrast with the no shear scenario in the second column, and the weak shear scenario in the third column. Each image is a view on the top of the domain and the magnitude of the total geostrophic velocity is being plotted. Simulation parameters: $N=64000$, $\eta=10^{-3}$, and $h=30$ min.}
    \label{fig:ShearEffect}
\end{figure}

Our formulation overcomes this critical limitation, allowing us to maintain accuracy and continue the simulation beyond the breakdown observed in their models. 
This advantage enables us to explore the extended dynamics of shearing effects on cyclogenesis over a 25-day period. 
In Figure~\ref{fig:ShearEffect}, we show how different initial shearing wind conditions ($A=-0.5$, $A=0$, and $A=0.1$) influence the intensity and organization of rotational systems. 
Our results reveal that with strong anticyclonic shear ($A=-0.5$), the formation of coherent rotational systems is significantly disrupted, while weaker or cyclonic shear ($A=0.1$) promotes the intensification of cyclonic structures, mirroring the findings of \cite{Wernli1998TheSimulations}, who observed pronounced differences in cyclone development based on the sign and magnitude of the imposed shear.

\begin{figure}[!ht]
    \centering
    \begin{subfigure}[b]{0.31\textwidth}
        \includegraphics[width=\textwidth]{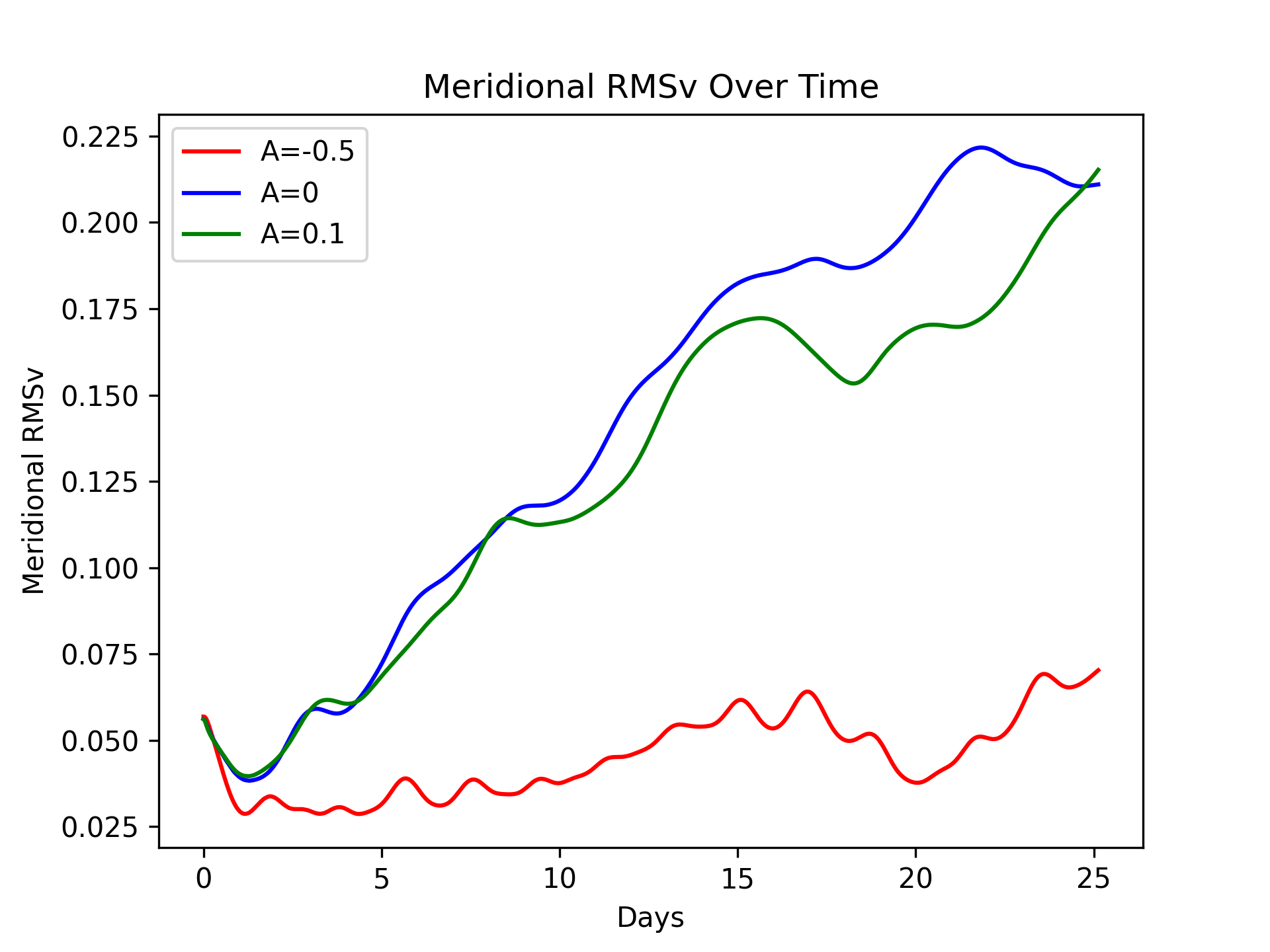}
        \caption{Meridional RMSv}
        \label{fig:MRMSv}
    \end{subfigure}
    \hfill 
    \begin{subfigure}[b]{0.31\textwidth}
        \includegraphics[width=\textwidth]{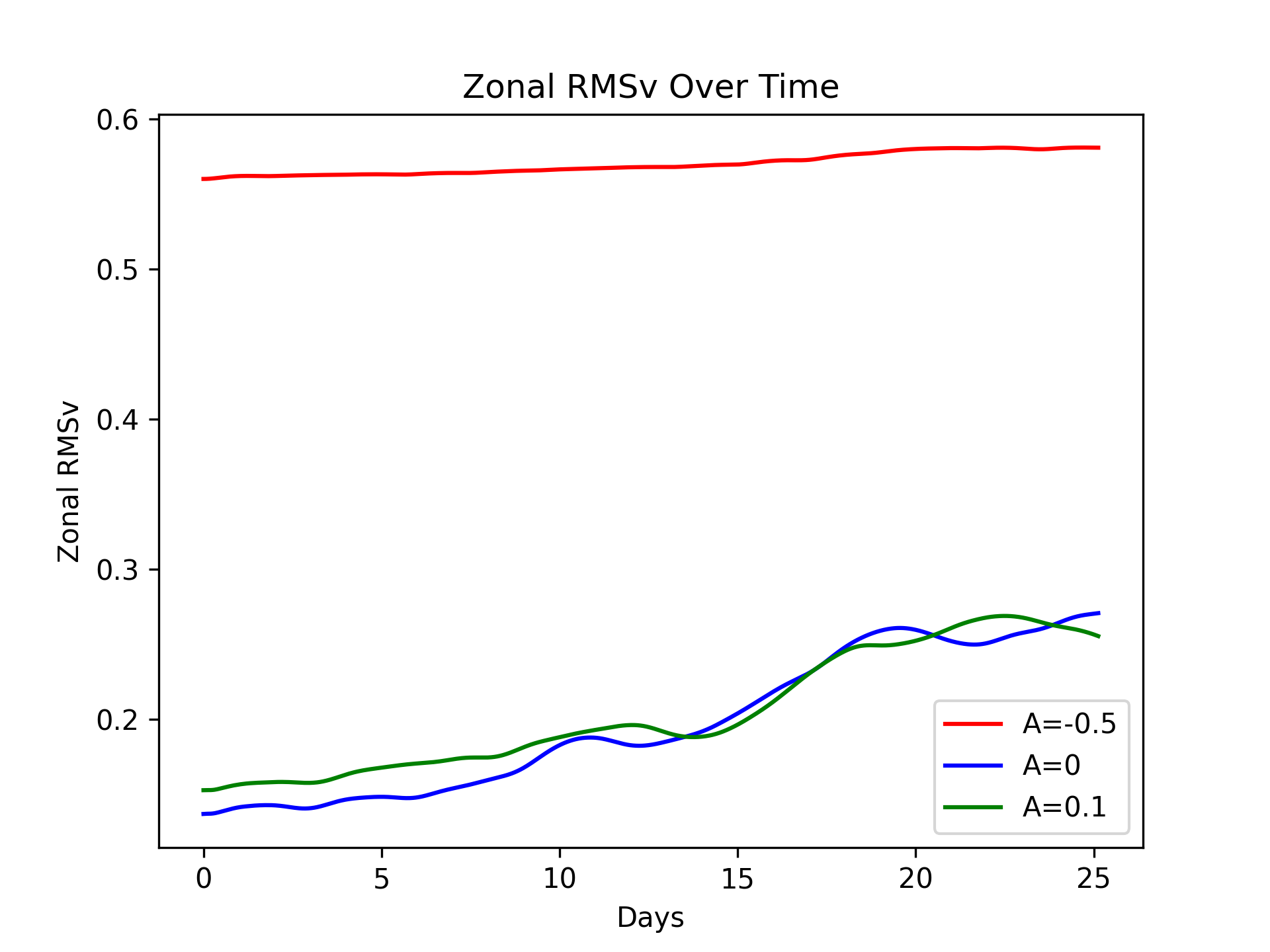}
        \caption{Zonal RMSv}
        \label{fig:ZRMSv}
    \end{subfigure}
    \hfill
    \begin{subfigure}[b]{0.31\textwidth}
        \includegraphics[width=\textwidth]{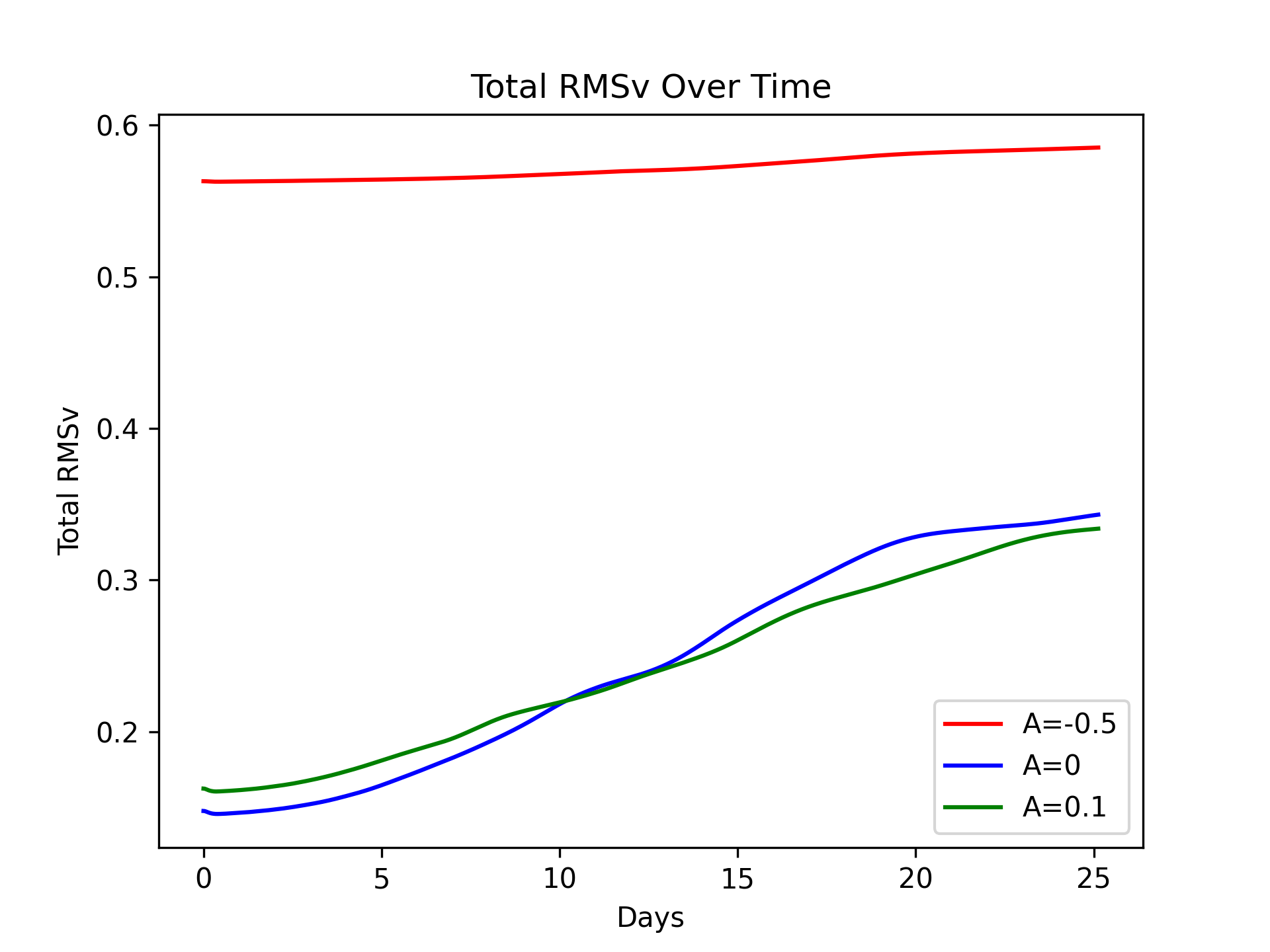}
        \caption{Total RMSv}
        \label{fig:TRMSv}
    \end{subfigure}
    \caption{Comparison of the evolution of root mean square velocity (RMSv) components under different initial shearing conditions over 25 days. \ref{fig:MRMSv} Meridional RMSv: shows the evolution of the y-direction velocity, with higher values indicating stronger meridional flow under weaker initial shear ($A = 0$, $A = 0.1$). \ref{fig:ZRMSv} Zonal RMSv: depicts the evolution of the x-direction velocity, demonstrating dominant zonal flow under strong initial shear ($A = -0.5$). \ref{fig:TRMSv} Total RMSv: presents the combined effect of both meridional and zonal components, illustrating the overall system dynamics under varying shearing conditions.}
    \label{fig:ShearWinds}
\end{figure}

Furthermore, as shown in Figure~\ref{fig:ShearWinds}, which separates the wind velocity into zonal and meridional components, the RMSv analysis highlights how shear conditions influence the balance between these components. 
For strong initial shear ($A=-0.5$), our results confirm that the zonal velocity dominates, leading to less organized and smaller-scale rotational systems, consistent with the findings by \cite{Davies1991TheDevelopment} where anticyclonic shear favored elongated cold fronts and weaker cyclones. 
Conversely, when shear is weak or absent ($A=0$ and $A=0.1$), the meridional velocity gains prominence, enhancing the development of more coherent cyclonic and anticyclonic structures, in line with the cyclonic shear experiments reported by \cite{Wernli1998TheSimulations}.

The ability to extend our simulations well beyond the timescales considered in previous studies provides new insights into the stability and evolution of these systems under sustained shear conditions. 
This prolonged analysis underscores the critical role of initial shearing in dictating the long-term behaviour of cyclonic structures, offering valuable extensions to the meteorological applications highlighted in \cite{Davies1991TheDevelopment} and \cite{Wernli1998TheSimulations}. 
Our results not only corroborate the sensitivity to shear observed in these foundational studies but also extend the understanding of how these dynamics evolve over longer timescales, providing a richer perspective on the impact of horizontal shear in geophysical fluid systems.

While providing detailed interpretations of the meteorological significance of specific atmospheric structures goes beyond the scope of the present work, our results clearly demonstrate the utility of our optimal transport formulation. 
By replicating the key features of the simulations by \cite{Wernli1998TheSimulations} and \cite{Davies1991TheDevelopment}, we highlight the potential of our approach as a valuable tool for investigating a broad range of SG atmospheric phenomena over extended timescales, enabling the study of complex dynamics that were previously inaccessible due to the limitations of traditional numerical formulations.
Indeed, we have adopted these initial conditions not to assert their definitive physical validity over indefinitely long timescales, but rather to serve as a qualitative benchmark for comparing SG solvers in three dimensions, a regime where, to our knowledge, few established benchmarks exist.

In standard Eulerian simulations, it is common to explicitly diagnose the ageostrophic velocity, $\vb{u}_{ag}$, to monitor the validity of the SG approximation (i.e., verifying $\norm{\vb{u}_{ag}}\ll\norm{\vb{u}_g}$). 
In the present semi-discrete optimal transport framework, however, extracting such a diagnostic is non-trivial and the reconstruction of the 3D physical velocity field  is the object of current research, and we cannot yet provide a rigorous estimate for $\vb{u}_{ag}$

Note however that the physical validity of the simulation is supported by the method's structural properties. 
The scheme maintains robust energy conservation throughout the 25-day evolution, and the resulting flow structures, specifically the formation and propagation of cyclones, qualitatively match the theoretical expectations and shorter-time Eulerian benchmarks established in the literature. 

\section*{Acknowledgements} 
We thank Jacques Vanneste and Mike Cullen for stimulating discussions related to this work and Quentin M\'erigot and Hugo Leclerc for supporting our use of their code, \textit{pysdot} (https://github.com/sd-ot/pysdot), and their continuous updates. 
We also want to give a special thanks to David Bourne, Beatrice Pelloni, and Charlie Egan, for their advice and guidance. 
TPL is supported by The Maxwell Institute Graduate School in Modelling, Analysis, and Computation, a Centre for Doctoral Training funded by the EPSRC (grant EP/S023291/1), the Scottish Funding Council, Heriot-Watt University and the University of Edinburgh.

\bibliographystyle{abbrv}
\bibliography{references}

\newpage
\appendix

\section{The Monge-Amp\`ere and Laplace Equations}\label{sec:MA-LP}

In this section, we justify the use of Laplace's equation for establishing our initial condition. 
The Monge-Ampère equation, related to the optimal transport map for the quadratic cost \cite{Brenier1991PolarFunctions}, plays a central role in coupling the SG equations with a transport equation. 
If $\nabla u_\# f = g$, then the potential $u$ satisfies the Monge-Amp\`ere equation given by
\begin{equation}\label{eq:MAE}
    f(x) = g(\nabla u(x)) \mathrm{det} D^2 u(x),
\end{equation}
for source probability measure $f:\R^n\to\R$, target probability measure $g:\R^n\to\R$, and a convex function $u\in\C^2(\X)$. 
Linearizing the Monge-Ampère equation around the quadratic potential leads to Poisson's equation, where the right-hand side depends on the gradient of the target measure. 
By neglecting this right-hand side, we obtain Laplace's equation, which is often used to approximate the initial condition in the quasi-geostrophic approximation.
In the case of an incompressible fluid the source measure is the Lebesgue measure so $f(\vb{x})=1$. We linearise Eq.~\eqref{eq:MAE} about 
\begin{align*}
    u(\vb{x})=\frac{1}{2}\vb{x}^T\vb{x},
\end{align*}
by adding the small perturbation $\varepsilon\Phi(\vb{x})$ to $u(\vb{x})$ to get
\begin{align*}
    1 &= g\qty(\vb{x}+\varepsilon\nabla\Phi(\vb{x}))\det\qty(\vb{I}+\varepsilon D^2\Phi(\vb{x})),
\end{align*}
and then differentiating both sides with respect to $\varepsilon$ to get
\begin{align*}
    0 &= \det\qty(\vb{I}+\varepsilon D^2\Phi(\vb{x}))\dv{}{\varepsilon}g\qty(\vb{x}+\varepsilon\nabla\Phi(\vb{x}))+g\qty(\vb{x}+\varepsilon\nabla\Phi(\vb{x}))\dv{}{\varepsilon}\det\qty(\vb{I}+\varepsilon D^2\Phi(\vb{x})) \\
    &= \det\qty(\vb{I}+\varepsilon D^2\Phi(\vb{x}))\nabla g(\vb{x}+\varepsilon\nabla\Phi(\vb{x}))\cdot\nabla\Phi(\vb{x}) \\
    &\quad\quad+g(\vb{x}+\varepsilon\nabla\Phi(\vb{x}))\det\qty(\vb{I}+\varepsilon D^2\Phi(\vb{x}))\tr\qty(\qty(\vb{I}+\varepsilon D^2\Phi(\vb{x}))^{-1}\cdot D^2\Phi(\vb{x})).
\end{align*}
By setting $\varepsilon =0$  and rearranging we arrive at
\begin{align*}
    \Delta\Phi(\vb{x}) = -\frac{1}{g(\vb{x})}\nabla g(\vb{x})\cdot\nabla\Phi(\vb{x}).
\end{align*}
Now by neglecting the right hand side we have
\begin{align*}
    \Delta \Phi(\vb{x}) = 0.
\end{align*}

This resulting Laplace equation is solved to derive the initial condition for the isolated large-scale tropical cyclone following the lead of \cite{Schar:1993}.

\section{Explicit solution for the perturbation} \label{sec:PerturbSoln}

In this section, we solve Laplace’s equation for the modified pressure, $\Phi$, decomposed as 
\begin{align*}
    \Phi(x_1, x_2, x_3)=\overline{\Phi}(x_1, x_2, x_3)+\widetilde{\Phi}(x_1, x_2, x_3),
\end{align*}
where $\overline{\Phi}$ is the background or steady state modified pressure and $\widetilde{\Phi}$ is the perturbed modified pressure. 
We do this to propagate the perturbation on the surfaces through the bulk of the domain. 
We consider a cuboid domain, subject to periodic boundary conditions in two directions and Neumann boundary conditions in the third direction. 
We also ensure that the compatibility condition for the Neumann problem is satisfied before proceeding with the solution. We begin with the Laplace equation
\begin{align*}
    \Delta \Phi(\vb{x}) = 0,
\end{align*}
in the cuboid domain $[-a,a]\times[-b,b]\times[0,c]$, with periodic boundary conditions in the $x_1$ and $x_2$ directions, and Neumann boundary conditions in the $x_3$ direction
\begin{align*}
    \pdv{\widetilde{\Phi}}{x_3}\eval_{x_3=0}=0.15h(x_1,x_2)\qq{and}\pdv{\widetilde{\Phi}}{x_3}\eval_{x_3=c}=-0.6h(x_1+1,x_2),
\end{align*}
where $h(x_1, x_2)$ is given by Eq.~\eqref{eq:Perturbation}. 
Note that $\overline{\Phi}$, introduced in Eq.~\eqref{eq:SteadyState}, is harmonic. 
Thus, we only need to solve
\begin{align*}
    \Delta \widetilde{\Phi}(\vb{x}) = 0.
\end{align*}
Before solving Laplace’s equation for $\widetilde\Phi$, we need to ensure the compatibility condition is satisfied. 
The compatibility condition requires that
\begin{align*}
    I_0+I_c=0,
\end{align*}
where
\begin{align*}
    I_0&=\int_{x_1=-a}^a\int_{x_2=-b}^b0.15h(x_1, x_2) \,\dd x_2 \dd x_1, \\
    I_c&=-\int_{x_1=-a}^a\int_{x_2=-b}^b0.6h(x_1 + 1, x_2) \,\dd x_2 \dd x_1.
\end{align*}
This condition is necessary for the solvability of the Neumann problem. 
However, the Neumann boundary conditions in \cite{Schar:1993} do not satisfy this condition. 
Therefore, we adjust the boundary conditions to: 
\begin{align*}
    \pdv{\widetilde{\Phi}}{x_3}\eval_{x_3=0}=0.15h(x_1,x_2)-\frac{I_0}{4ab}\qq{and}\pdv{\widetilde{\Phi}}{x_3}\eval_{x_3=c}=-0.6h(x_1+1,x_2)-\frac{I_c}{4ab}.
\end{align*}
With the compatibility condition now satisfied, we can proceed to solve Laplace's equation.
We start by making the usual ansatz and expanding $\widetilde{\Phi}(x_1,x_2,x_3)$ in a Fourier series
\begin{align*}
    \widetilde{\Phi}(x_1,x_2,x_3) = \sum_{n=-\infty}^{\infty} \sum_{m=-\infty}^{\infty} \exp \left( \frac{\pi i n x_1}{a} \right) \exp \left( \frac{\pi i m x_2}{b} \right) Z_{n,m}(x_3). 
\end{align*}
where $Z_{n,m}(x_3)$ are the unknown coefficient functions to be determined. 
By substituting this ansatz into the Laplace equation $\Delta \widetilde{\Phi}(\vb{x})=0$, we obtain the following set of ordinary differential equations for $Z_{n,m}(x_3)$:
\begin{align*}
    \dv[2]{}{x_3}Z_{n,m}(x_3)=k^2_{n,m}Z_{n,m}(x_3),
\end{align*}
where $k_{n,m}=\pi\sqrt{\qty(\frac{n}{a})^2+\qty(\frac{m}{b})^2}$. 
The general solution for $Z_{n,m}(x_3)$ is
\begin{align*}
    Z_{n,m}(x_3)=C_{n,m}\exp\qty(k_{n,m}x_3)+D_{n,m}\exp\qty(-k_{n,m}x_3),
\end{align*}
where $C_{n,m}$ and $D_{n,m}$ are constants determined by the boundary conditions. 
Thus the solution is
\begin{align*}
    \widetilde{\Phi}(x_1,x_2,x_3)&=\sum_{n=-\infty}^{\infty}\sum_{m=-\infty}^{\infty}\exp\qty(\frac{\pi i n x_1}{a})\exp\qty(\frac{\pi i m x_2}{b}) \qty(C_{n,m}\exp(k_{n,m}x_3)+D_{n,m}\exp(-k_{n,m}x_3)).
\end{align*}
Next, we apply the Neumann boundary conditions at $x_3=0$ and $x_3=c$ to find $C_{n,m}$ and $D_{n,m}$. 
First, we compute the derivative of $\widetilde{\Phi}$ with respect to $x_3$
\begin{align*}
    \pdv{\widetilde{\Phi}}{x_3}(x_1,x_2,x_3)=\sum_{n=-\infty}^{\infty}\sum_{m=-\infty}^{\infty}k_{n,m}\exp\qty(\frac{\pi i n x_1}{a})\exp\qty(\frac{\pi i m x_2}{b}) \qty(C_{n,m}\exp(k_{n,m}x_3)-D_{n,m}\exp(-k_{n,m}x_3)).
\end{align*}
Now, we find the Fourier transforms of the boundary conditions. 
For the lower boundary $x_3=0$, the Fourier coefficients $A_{n,m}$ are given by
\begin{align*}
    A_{n,m}=\frac{1}{4 ab}\int_{x_1=-a}^a\int_{x_2=-b}^b\qty(0.15h(x_1,x_2)-\frac{I_0}{4ab})\exp\qty(-\frac{\pi i n x_1}{a})\exp\qty(-\frac{\pi i m x_2}{b})\,\dd x_2\dd x_1.
\end{align*}
For the upper boundary $x_3=c$, the Fourier coefficients $B_{nm}$ are given by
\begin{align*}
    B_{n,m}&=\frac{1}{4 ab}\int_{x_1=-a}^a\int_{x_2=-b}^b \qty(-0.6h(x_1+1,x_2)-\frac{I_c}{4ab})\exp\qty(-\frac{\pi i n x_1}{a})\exp\qty(-\frac{\pi i m x_2}{b})\,\dd x_2\dd x_1.
\end{align*}
We then solve the following system of equations to determine $C_{nm}$ and $D_{nm}$ :
\begin{align*}
    A_{nm}&=k_{nm}\qty(C_{nm}-D_{nm}), \\
    B_{nm}&=k_{nm}\qty(C_{nm}\exp\qty(k_{nm}c)-D_{nm}\exp\qty(-k_{nm}c)).
\end{align*}
For the case $n=m=0$, we set $C_{00}=D_{00}=0$, which corresponds to the average perturbation on the surfaces of the domain. 
With the coefficients $C_{nm}$ and $D_{nm}$ determined, the initial condition is fully established, and we are now ready to proceed with the numerical solution of the system.

\end{document}